\input ./style/arxiv-general.cfg
\documentclass[aop,MSNbibl,dvips]{arximspdf}
\makeatletter
   \@ifpackageloaded{graphicx}{}{\usepackage{graphicx}}
\makeatother

\doi{10.1214/14-AOP929} 
\volume{43}
\issue{4}
\pubyear{2015}
\firstpage{2119}
\lastpage{2150}

\makeatletter

\newcommand{\rright}{\right}
\newcommand{\lleft}{\left}
\newcommand{\rrVert}{\Vert}
\newcommand{\rrvert}{\vert}
\newcommand{\llVert}{\Vert}
\newcommand{\llvert}{\vert}
\newtheorem{theorem}{Theorem}[section]
\newtheorem{lemma}[theorem]{Lemma}
\newtheorem{proposition}[theorem]{Proposition}

\newproclaim{ex}{Example}
\newproclaim{rem}{Remark}
\makeatother

\begin{document}
\begin{frontmatter}

\title{Subordination for the sum of two random matrices}
\runtitle{Subordination for the sum of two random matrices}

\begin{aug}
\author[A]{\fnms{V.}~\snm{Kargin}\corref{}\ead[label=e1]{vladislav.kargin@gmail.com}}
\runauthor{V. Kargin}
\affiliation{University of Cambridge}
\address[A]{Statistical Laboratory\\
Center for Mathematical Sciences\\
Wilberforce Road\\
Cambridge CB3 OWB\\
United Kingdom\\
\printead{e1}} 
\end{aug}

\received{\smonth{3} \syear{2013}}
\revised{\smonth{1} \syear{2014}}

%
\begin{abstract}
This paper is about the relation of random matrix theory and the
subordination phenomenon in complex analysis. We find that the
resolvent of
the sum of two random matrices is approximately subordinated to the
resolvents of the original matrices. We estimate the error terms in this
relation and in the subordination relation for the traces of the
resolvents. This allows us
to prove a local limit law for eigenvalues and a delocalization result for
eigenvectors of the sum of two random matrices. In addition, we use
subordination to determine the limit of the largest eigenvalue for the
rank-one deformations of unitary-invariant random
matrices.
\end{abstract}

\end{frontmatter}

\section{Introduction}\label{sec1}

\subsection{Subordination}\label{sec11}

Much of the modern approach to random matrices is based on the analysis of
how the resolvent of a matrix $A$, that is, the function $G_{A} (
z ) = ( A-zI ) ^{-1}$, behaves when $A$ is modified by a
random perturbation (see~\cite{erdosschleinyau09}, e.g.). In
this paper, we investigate what happens with the resolvent if an independent
rotationally invariant random matrix $B$ is added to $A$. We find that the
resolvent of the sum $A+B$ is (approximately) subordinated to the resolvent
of the original matrix $A$.

The concept of subordination comes from the complex analysis. If $f(z)$
and $%
g(z)$ are two functions which are analytic in the upper half-plane
$\mathbb{C%
}^{+}:=\{z\dvtx \operatorname{Im}z>0\}$, then $f(z)$ is \emph{subordinated} to $g(z)$
if there exists an analytic function $\omega(z)\dvtx \mathbb
{C}^{+}\rightarrow
\mathbb{C}^{+}$, such that $f(z)=g(\omega(z))$ and $\operatorname{Im}\omega
(z)\geq\operatorname{Im}z$ for all $z\in\mathbb{C}^{+}$. In this
definition, $%
f(z)$ and $g(z)$ can be vector or operator valued functions.

Voiculescu and Biane \cite{voiculescu93,biane98b} have discovered
that the subordination holds for the resolvent of the sum of two free
operators in a von Neumann algebra. (See also \cite{belinschibercovici07}
and \cite{chistyakovgotze11} for different proofs of these results.) This
subordination result can be formulated as follows (cf. Theorem 3.1 in
\cite{biane98b}). Let $\mathcal{A}$ be a von Neumann operator algebra with the
normal faithful trace $\tau\dvtx \mathcal{A}\rightarrow\mathbb{C}$. If two
self-adjoint operators $A,B\in\mathcal{A}$ are free in the sense of
Voiculescu (see \cite{nicaspeicher06}), then the following identity holds:
%
\begin{equation}
\tau\bigl(G_{A+B} ( z ) |A\bigr)=G_{A}\bigl(
\omega_{B} ( z ) \bigr), \label{subordinationoperator}
\end{equation}
where $\tau(\cdot|A)$ denotes the conditional expectation on the
subalgebra generated by operator $A$, and $\omega_{B}(z)$ is a function
analytic in $\mathbb{C}^{+}$ and such that $\operatorname{Im}\omega
_{B}(z)\geq
\operatorname{Im}z$. In other words, $\tau(G_{A+B}(z)|A)$ is subordinated
to $%
G_{A}(z)$.

This subordination result is very useful since it implies results about the
smoothness of the spectral distribution of the sum $A+B$.

Since large independent random matrices are asymptotically free \cite{voiculescu91,speicher93}, it is natural to ask whether
subordination holds in the context of random matrices. Some results in this
direction have been recently obtained in \cite{capitainedonati-martinferalfevrier11,male11} and \cite{capitaine12}. In \cite{capitainedonati-martinferalfevrier11} (which
builds on an earlier work in \cite{capitainedonati-martinferal09}), the
authors study the matrix $A_{N}+W_{N}/\sqrt{N}$, where $A_{N}$ and
$W_{N}$ are $N$-by-$N$ Hermitian matrices, $A_{N}$ is deterministic
and $W_{N}$
is Wigner. It is assumed that the eigenvalue distribution of $A_{N}$ weakly
converges to a measure $\nu$ as $N\rightarrow\infty$ and that the largest
$r$ eigenvalues of $A_{N}$ (``spikes'') are
fixed and are outside of the support of $\nu$. The authors are interested
in the behavior of $r$ largest eigenvalues of $A_{N}+W_{N}/\sqrt{N}$ and
this question leads them to the study of the subordination for the
trace of
the resolvent of $A_{N}+W_{N}/\sqrt{N}$.

In further developments, in \cite{male11} and in \cite{capitaine12}, the
setup of \cite{capitainedonati-martinferalfevrier11} is generalized for
perturbations of the block random matrices and sample covariance matrices,
respectively.

We are interested in a somewhat different setup. Let $\widetilde{A}$
and $%
\widetilde{B}$ be two \mbox{$N$-by-$N$} diagonal matrices with real entries. Define
the random matrices $A:=V\widetilde{A}V^{\ast}$ and $B:=U\widetilde
{B}%
U^{\ast}$ where $U$ and $V$ are two $N$-by-$N$ random independent uniformly
distributed unitary matrices, and define $H:=A+B$. Note that the
distribution of eigenvalues of $H$ is the same as that of $\widetilde
{A}+U%
\widetilde{B}U^{\ast}$, however, it will be convenient to treat $A$
and $B$
symmetrically. The \emph{resolvent} of $H$ is defined as $G_{H} (
z ):= ( H-zI ) ^{-1}$ and the \emph{Stieltjes
transform} of $%
H $ is defined as the normalized trace of the resolvent:
\[
m_{H} ( z ):=\frac{1}{N}\operatorname{Tr}\bigl(G_{H} ( z
) \bigr).
\]
The resolvents and the Stieltjes transforms of matrices $A$ and $B$ are
defined similarly.

Is it true that $G_{H} ( z ) $ is subordinated to
$G_{A} (
z ) $ and $G_{B} ( z ) $ for sufficiently large~$N$?

First, we need to define a candidate subordination function. Let
%
\begin{equation}
\omega_{B} ( z ):=z-\frac{\mathbb{E}f_{B} ( z
) }{%
\mathbb{E}m_{H} ( z ) }\quad\mbox{and}\quad\omega
_{A} ( z ):=z-%
\frac{\mathbb{E}f_{A} ( z ) }{\mathbb{E}m_{H} (
z ) }, \label{omegadefinition}
\end{equation}
where
\[
f_{B} ( z ):=N^{-1}\operatorname{Tr} \bigl( BG_{H} ( z
) \bigr)\quad\mbox{and}\quad f_{A} ( z ):=N^{-1}\operatorname{Tr}
\bigl( AG_{H} ( z ) \bigr).
\]

We claim $\omega_{A} ( z ) $ and $\omega_{B} (
z ) $ are
``almost'' subordination functions.

\begin{theorem}
\label{Propsubordination} Assume that $\eta:=\operatorname{Im}z\in (
0,1 ) $ and $\llvert \operatorname{Re}z\rrvert \leq K (
A,B ):=\max \{ \llVert  A\rrVert,\llVert  B\rrVert
 \}. $
Then for all $N\gg\eta^{-5}$,
\[
\min \bigl\{ \operatorname{Im} \bigl( \omega_{A} ( z ) \bigr),
\operatorname{Im} \bigl( \omega_{B} ( z ) \bigr) \bigr\} \geq
\eta-\frac{c}{N\eta^{7}},
\]
with $c>0$ that depends only on $K ( A,B ) $.
\end{theorem}

In other words, for all sufficiently large $N$, the excess of the imaginary
parts of functions $\omega_{A} ( z ) $ and $\omega
_{B} (
z ) $ over $\operatorname{Im}z$ is almost nonnegative. The proof of this
theorem is postponed to the next section.

Now we are able to formulate the main result.

[We use the following notation. The average of a random variable $X$
over $%
U $ is denoted by $\mathbb{E}_{U} ( X ):=\mathbb{E} (
X|V ) $ and the average over $V$ is denoted by $\mathbb
{E}_{V} (
X ):=\mathbb{E} ( X|U ) $. The unconditional expectation
value is denoted by $\mathbb{E} ( X ) $. Similar notation
will be
used for conditional probabilities and variances. For example, $\mathbb
{V}%
\mathrm{ar}_{V} ( X ) =\mathbb{E}_{V} (  (
X-\mathbb{E}%
_{V}X ) ^{2} ) $. The notation $x\ll y$ and $x=O (
y ) $
mean that there exists a constant $C>0$ such that $\llVert  x\rrVert
\leq Cy$. The constants in these inequalities may depend on $K (
A,B ):=\max \{ \llVert  A\rrVert,\llVert
B\rrVert
 \} $. The norm $\llVert \cdot\rrVert $ is the usual uniform
norm on matrices.]

\begin{theorem}
\label{PropRAestimate} Assume that $\eta:=\operatorname{Im}z\in (
0,1 ) $ and $\llvert \operatorname{Re}z\rrvert \leq K (
A,B ):=\max \{ \llVert  A\rrVert,\llVert  B\rrVert
 \} $.
Suppose that $N\gg\eta^{-7}$. Then we have:
\begin{longlist}[(iii)]
\item[(i)]
%
\begin{equation}
\mathbb{E}_{U}G_{H} ( z ) -G_{A} \bigl(
\omega_{B} ( z ) \bigr) =O \biggl( \frac{1}{N\eta^{6}} \biggr),
\end{equation}

\item[(ii)]
\[
\mathbb{E}m_{H} ( z ) -m_{A} \bigl( \omega_{B}
( z ) \bigr) =O \biggl( \frac{1}{N^{2}\eta^{6}} \biggr),
\]

\item[(iii)]
\[
\mathbb{P}_{U} \bigl\{ \bigl\llvert \bigl( G_{H} ( z )
\bigr) _{ij}- \bigl( G_{A} \bigl( \omega_{B} ( z )
\bigr) \bigr) _{ij}\bigr\rrvert \geq\delta \bigr\} \leq\exp \bigl( -c
\delta ^{2}N^{3/7} \bigr),
\]
for all $N\geq N_{0}$, where $N_{0}$ can depend on $K ( A,B
) $ and
on $\delta$.
\end{longlist}
\end{theorem}

Estimates similar to estimates in parts (i) and (iii) hold for the
conditional expectations with respect to $V$. The expectation in part (ii)
is unconditional since $\mathbb{E}m_{H} ( z ) =\mathbb{E}%
_{U}m_{H} ( z ) =\mathbb{E}_{V}m_{H} ( z ) $ (by an
application of Lemma \ref{lemmaexpectationandtrace} in the
\hyperref[app]{Appendix}). An
estimate similar to the estimate in part (ii) holds for $\mathbb{E}%
m_{H} ( z ) -m_{B} ( \omega_{A} ( z )
 ) $.

It is easy to check (see Lemma \ref{PropEGdiagonal}) that if the
basis is
chosen in such a way that $A$ is diagonal then the matrix $\mathbb{E}%
_{U}G_{H} ( z ) $ is also diagonal. Hence, parts (i) and
(iii) of
the theorem essentially say that if $A$ is diagonal then $G_{H} (
z ) $ is approximately diagonal and its diagonal entries satisfy the
formula
\[
\bigl( G_{H} ( z ) \bigr) _{ii}\approx\frac
{1}{A_{ii}-\omega
_{B} ( z ) }.
\]

Part (ii) says that taking the trace of the resolvent makes the error
in the
approximate formula even smaller, $O ( N^{-2}\eta^{-6} ) $ instead
of $O ( N^{-1}\eta^{-6} ) $.

It is interesting to compare this result with the results in \cite{capitainedonati-martinferalfevrier11}. Let $H=A+W/\sqrt{N}$ where
$W$ is
the $N$-by-$N$ Wigner matrix with i.i.d. Gaussian entries of variance $%
\sigma^{2}$. Let $m_{H} ( z ):=N^{-1}\mathbb{E}\operatorname{Tr}%
G_{H} ( z ) $. (Note that $\mathbb{E}$ has a slightly different
meaning here. It is the expectation taken with respect to the
randomness in
the Wigner matrix.) It was proved in \cite{capitainedonati-martinferalfevrier11} that for all $z\in\mathbb{C}^{+}$,
%
\begin{equation}
m_{H} ( z ) =m_{A} \bigl( z+\sigma^{2}m_{H}
( z ) \bigr) +O \bigl( N^{-2} \bigr), \label{formulasubordinationcapitaine}
\end{equation}
with the constant in the $O$-term that depends on $z$. In addition, if $W$
is a non-Gaussian Wigner matrix, then the same formula is proved in
\cite{capitainedonati-martinferalfevrier11} with an additional term on the
right of the form $L ( z ) /N$.

Formula (\ref{formulasubordinationcapitaine}) gives a subordination result
with the subordination function $\omega_{B}=z+\sigma^{2}m_{H} (
z ) $. Our Theorem \ref{PropRAestimate} holds for a more general
matrix model and gives estimates for resolvents as well as for the Stieltjes
transforms. These estimates give the explicit dependence of the error term
on $z$ unlike formula (\ref{formulasubordinationcapitaine}). These
advantages are crucial for the applications to the local distribution of
eigenvalues and delocalization of eigenvectors.

Now we turn to these applications.

\subsection{Delocalization}\label{sec12}

\emph{Delocalization} of eigenvectors generally refers to the
situation when
all individual coordinates of a normalized eigenvector $v_{\alpha}$ in a
specific basis are not greater than $N^{-\kappa/2}$ with a high
probability. Because of normalization, the eigenvector is forced to be
spread over at least $N^{\kappa}$ coordinates and it is customary to say
that the delocalization length of eigenvectors is at least $N^{\kappa}$.

The question about delocalization of eigenvectors frequently occurs in
physics. For example, a famous open problem is to show that for $d\geq3$
the eigenvectors of random Schr\"{o}dinger operators on $\mathbb
{Z}^{d}$ are
completely delocalized for small disorder. Recently, there was some progress
on delocalization of eigenvectors in simpler models, for instance, in the
case of random Wigner matrices and in the case of random band matrices. In
the former case, the complete delocalization has been established recently
(see \cite{erdos10} for a review) and the method is similar to the method
that is used in this paper. In the case of band matrices, it is expected
that complete delocalization holds for matrices with the band width $W$
greater than $\sqrt{N}$. What was actually shown in this case is that the
delocalization length is greater than $W^{1+d/6}$ (\cite{erdosknowles11},
and an improvement was recently achieved in \cite{ekyy12}). The method is
based on quantum diffusion and different from the method that is used in
this paper.

In our model, we say that the eigenvectors $\nu_{a}^{(N)}$ of a
sequence of
matrices $H_{N}=A_{N}+UB_{N}U^{\ast}$ are \emph{delocalized} at
length $%
N^{\kappa}$ in the interval $I$, if there exists $\delta>0$ such that
\[
\mathbb{P} \bigl\{ \bigl\llvert v_{a}^{ ( N ) } ( i ) \bigr
\rrvert ^{2}>N^{-\kappa}\log N \bigr\} \leq\exp \bigl(
-N^{-\delta
} \bigr),
\]
for all sufficiently large $N$, all $i\in \{ 1,\ldots,N \}
$ and
all $\nu_{a}^{(N)}$ such that the corresponding eigenvalues are in the
interval $I$.

Let $\mu_{A_{N}}$ be the empirical measure of eigenvalues of $A_{N}$, that
is, $\mu_{A_{N}}:=N^{-1}\sum_{k=1}^{N}\delta_{\lambda_{k}}$, where $
\lambda_{k}$ are eigenvalues of $\mu_{A_{N}}$. Define $\mu_{B_{N}}$
similarly. We are going to prove that the eigenvalues of $H_{N}$ are
delocalized at a certain scale if $\mu_{A_{N}}$ and $\mu_{B_{N}}$ are
close enough to a couple of measures that satisfy a regularity conditions.

As a measure of closeness between probability measures $\mu$ and $\nu
$, we
use the \emph{L\'{e}vy distance}
\[
d_{L} ( \mu,\nu ) =\sup_{x}\inf \bigl\{ s
\geq0\dvtx F_{\nu
} ( x-s ) -s\leq F_{\mu} ( x ) \leq
F_{\nu} ( x+s ) +s \bigr\},
\]
where\vspace*{2pt} $F_{\mu} ( t ) $ and $F_{\nu} ( t ) $ are the
cumulative distribution functions of $\mu$ and $\nu$. Note that $\mu
^{ ( N ) }\rightarrow\mu$ in distribution if and only if $%
d_{L} ( \mu^{ ( N ) },\mu ) \rightarrow0$.
(See Theorem~III.1.2 on page~314 and Exercise III.1.4 on page~316 in \cite{shiryaev96}.)

\begin{theorem}
\label{Thmdelocalization}Assume that \textup{(i)} a pair of probability
measures $%
 ( \mu_{\alpha},\mu_{\beta} ) $ is smooth in a closed
interval $%
I$, and \textup{(ii)} for a sequence of $A_{N}$ and $B_{N}$, $\max \{
\llVert
A_{N}\rrVert,\llVert  B_{N}\rrVert  \} \leq K$ for
all $N$.
Then there exists $s>0$ such that if
%
\begin{equation}
\max \bigl\{ d_{L}(\mu_{A_{N}},\mu_{\alpha}),d_{L}(
\mu_{B_{N}},\mu _{\beta}) \bigr\} \leq s \label{conditionclosenessAB}
\end{equation}
for\vspace*{1pt} $N$ large enough, then eigenvectors of $H_{N}=A_{N}+UB_{N}U^{\ast
}$ are
delocalized at scale $N^{1/7}$ in the interval $I$.
\end{theorem}

Note that we do not require the measures $\mu_{A_{N}}$ and~$\mu_{B_{N}}$
to converge to $\mu_{\alpha}$ and~$\mu_{\beta}$. It is enough that they
are sufficiently close to $\mu_{\alpha}$ and $\mu_{\beta}$ for all large~$N$. The reason for this is that this weaker condition is enough to ensure
that the subordination functions of the pair $ ( \mu_{A_{N}},\mu
_{B_{N}} ) $  are separated from zero in the region $\eta\geq
N^{-1/7}$. On the other hand, if $\mu_{A_{N}}$ and $\mu_{B_{N}}$ do
converge to $\mu_{\alpha}$ and~$\mu_{\beta}$, then condition
(\ref{conditionclosenessAB}) is automatically satisfied. This is perhaps the
most important case in applications.

We still need to explain what is meant by the smoothness of a pair
$ (
\mu_{\alpha},\mu_{\beta} ) $. Let $\mu_{\alpha}$ and $\mu
_{\beta
}$ be two\vspace*{2pt} probability measures with bounded support, and let $m_{\alpha
} ( z ):=\int ( t-z ) ^{-1}\mu_{\alpha} (
\,dt ) $ and $m_{\beta} ( z ):=\int ( t-z
) ^{-1}\mu
_{\beta} ( dt ) $. The system of equations
%
\begin{eqnarray}\label{system}
m ( z ) &=&m_{\alpha} \bigl( \omega_{\beta} ( z ) \bigr),\nonumber
\\
m ( z ) &=&m_{\beta} \bigl( \omega_{\alpha} ( z ) \bigr)\quad
\mbox{and}
\\
z-\frac{1}{m ( z ) } &=&\omega_{\alpha} ( z ) +\omega _{\beta} ( z
)\nonumber
\end{eqnarray}
has a unique solution $ ( m ( z ),\omega_{\alpha
} (
z ),\omega_{\beta} ( z )  ) $ in the class of
functions that are analytic in $\mathbb{C}^{+}= \{ z\dvtx \operatorname{Im}%
z>0 \} $ and that have the following expansions at infinity:
%
\begin{eqnarray}\label{asymptoticconditions}
m ( z ) &=&-z^{-1}+O \bigl( z^{-2} \bigr),
\nonumber\\[-8pt]\\[-8pt]
\omega_{\alpha} ( z ) &=&z+O ( 1 )\quad\mbox {and}\quad\omega
_{\beta} ( z ) =z+O ( 1 ).
\nonumber
\end{eqnarray}
The function $m ( z ) $, which we denote as $m_{\mu_{\alpha
}\boxplus\mu_{\beta}} ( z ) $, is the Stieltjes transform
of a
probability measure which is called the \emph{free convolution} of measures
$\mu_{\alpha}$ and $\mu_{\beta}$ and denoted $\mu_{\alpha
}\boxplus\mu
_{\beta}$. The functions $\omega_{\alpha} ( z ) $ and
$\omega
_{\beta} ( z ) $ are the subordination functions for the free
convolution.

By Theorem 3.3 in \cite{belinschi08}, the limits $\omega_{j}(x)=\lim_{\eta
\downarrow0}\operatorname{Im}\omega_{j} ( x+i\eta ) $ exist
for $%
j=\alpha,\beta$, and we make the following definition. A pair of
probability measures on the real line $ ( \mu_{\alpha},\mu
_{\beta
} ) $ is said to be \emph{smooth} at $x\in\mathbb{R}$ if the following
two conditions hold:
\begin{longlist}[(A)]
\item[(A)] $\operatorname{Im}\omega_{j}(x)>0 $ for $j=\alpha,\beta$, and

\item[(B)]
%
\begin{equation}\label{genericity}
k_{\mu}(x):=\frac{1}{m_{\mu_{\alpha}}^{\prime}(\omega_{\beta
}(x))}+%
\frac{1}{m_{\mu_{\beta}}^{\prime}(\omega_{\alpha}(x))}-\bigl(
\omega _{\alpha
}(x)+\omega_{\beta}(x)-x\bigr)^{2}\neq0.
\end{equation}
\end{longlist}

We say that the pair $ ( \mu_{\alpha},\mu_{\beta} ) $ is
\emph{smooth} \emph{in interval} $I\subset\mathbb{R}$ if $\omega_{\alpha
}(z) $ and $\omega_{\beta}(z)$ are continuous in a rectangle $ \{
z=x+i\eta
|x\in I,0\leq\eta\leq\varepsilon \} $ where $\varepsilon$ is a
positive constant, and if the pair $ ( \mu_{\alpha},\mu_{\beta
} ) $ is smooth at every point of $I$.

The proof of Theorem \ref{Thmdelocalization} is based on part
(iii) of
Theorem \ref{PropRAestimate} which imply that
$\operatorname{Im} ( G_{H} ) _{kk} ( \lambda_{a}+i\eta
 ) \leq
\operatorname{Im} ( G_{A} ( \omega_{B} ( \lambda_{a}+i\eta
 )
 )  ) _{kk}+\delta$. Then the assumption of smoothness leads
(after some work) to the conclusion that the quantity on the right is
bounded for all $k$ and all $N\gg\eta^{-1/7}$ with high probability.
Therefore, the components of the eigenvector corresponding to $\lambda_{a}$
can be estimated by using the~bound on the resolvent
%
\begin{equation}
\bigl\llvert v_{a} ( k ) \bigr\rrvert ^{2}\leq\eta\operatorname{Im}%
G_{kk} ( \lambda_{a}+i\eta ) \leq C\eta\leq
CN^{-1/7}\log N.
\end{equation}
To get the last inequality, $\eta$ is chosen as $N^{-1/7}\log N$ so that
Theorem \ref{PropRAestimate} is applicable. The details are
postponed to
Section~\ref{sec3}.

Let us add some comments about the assumption of smoothness. Condition~(B) is technical and holds for a generic point $x\in\mathbb{R}$. It ensures
that the solution of the system (\ref{system}) at $x$ is stable with respect
to a small perturbation in the system. Condition (A) is essential and
closely related to regularity properties of the measure $\mu_{\alpha
}\boxplus\mu_{\beta}$ at $x$. Here are some cases when it holds:
\begin{longlist}[(ii)]
\item[(i)] If $\mu_{\alpha}=\mu_{\beta}=\mu$, and $\mu\boxplus\mu$ is
absolutely continuous with positive density at $x$, then (A) is
satisfied at
the point $x$. In particular, if $\mu$ is an arbitrary measure that does
not have an atom with the mass greater than $1/2$, then condition
(A) is
satisfied at every point inside the support of $\mu\boxplus\mu$.

\item[(ii)] If one of the probability measures has the semicircle distribution with
the density $f_{\mathrm{sc}} ( x ) =\frac{1}{2\pi}\sqrt{ (
4-x^{2} ) _{+}}$, the density of $\mu_{\mathrm{sc}}\boxplus\mu$ is positive
at $x$, and $\llvert  m_{\mu_{\mathrm{sc}}\boxplus\mu} ( x )
\rrvert \neq1$, then condition (A) is satisfied. (For a more
detailed discussion of these examples, the reader can see Propositions 1.4
and 1.5 in \cite{kargin13a}.)
\end{longlist}

In fact, smoothness is likely to be a typical situation for pairs
$ (
\mu_{\alpha},\mu_{\beta} ) $. This is because the free convolution
operation has very strong smoothing properties. Even if we start with two
discrete measures $\mu_{\alpha}$ and $\mu_{\beta}$, the free convolution
$\mu_{\alpha}\boxplus\mu_{\beta}$ is absolutely continuous provided
that the masses of an atom of $\mu_{\alpha}$ and an atom of $\mu
_{\beta}$
do not add up to more than $1$.

How does one find pairs which are not smooth? A pair of measures $
( \mu
_{\alpha},\mu_{\beta} ) $ is not smooth at a point where the density
of $\mu_{\alpha}\boxplus\mu_{\beta}$ vanishes, in particular at the
boundary of the support of $\mu_{\alpha}\boxplus\mu_{\beta}$. One other
example occurs when both $\mu_{\alpha}$ and $\mu_{\beta}$ have an atom,
and the sum of the atoms' masses is greater than $1$. In this case, the free
convolution $\mu_{\alpha}\boxplus\mu_{\beta}$ also has an atom and,
therefore, the pair $ ( \mu_{\alpha},\mu_{\beta} ) $ is not
smooth at the location of this atom.

The result in Theorem \ref{Thmdelocalization} is certainly not
optimal. The
true localization length is probably of order $N$ under assumptions of the
theorem, that is, eigenvectors are likely to be completely delocalized.

\subsection{Local limit for eigenvalue distribution}\label{sec13}

Another consequence of Theorem~\ref{PropRAestimate} is the
convergence of
the eigenvalue counting measure on the local scale.

Let $\mathcal{N}_{\eta^{\ast}} ( x ) $ be the number of
eigenvalues of $H_{N}$ in the interval $I^{\ast}= [ x-\eta^{\ast
},x+\eta^{\ast} ] $. What can be said about $\mathcal{N}_{\eta
^{\ast
}} ( x ) / ( 2\eta^{\ast}N ) $ when
$N\rightarrow\infty
?$ If $\eta^{\ast}$ is fixed, then it is known \cite{voiculescu91}
and \cite{speicher93} that the limit approaches $\mu_{\alpha}\boxplus
\mu
_{\beta} ( I^{\ast} ) / ( 2\eta^{\ast} ) $. Local
limit theorems address the question of what happens if $\eta^{\ast}$ is
not fixed but approaches $0$ when $N\rightarrow\infty$.

\begin{theorem}
\label{Thmlocallaw}Assume that \textup{(i)} $\max \{ d_{L}(\mu
_{A_{N}},\mu
_{\alpha}),d_{L}(\mu_{B_{N}},\mu_{\beta}) \} \rightarrow0$,  \textup{(ii)}~the
pair of probability measures $ ( \mu_{\alpha},\mu_{\beta
} ) $ is smooth on interval $I$ and \textup{(iii)}~$\max \{ \llVert
A_{N}\rrVert,\llVert  B_{N}\rrVert  \} \leq K$ for
all $N$.
Let $\rho_{\mu_{\alpha}\boxplus\mu_{\beta}}$ denote the density of
$\mu_{\alpha}\boxplus\mu_{\beta}$, and let $\eta^{\ast
}=cN^{-1/7}\log
N $. Then, for every $x\in I$,
\[
\frac{\mathcal{N}_{\eta^{\ast}} ( x ) }{2\eta^{\ast}N}%
\rightarrow\rho_{\mu_{\alpha}\boxplus\mu_{\beta}} ( x )
\]
in probability.
\end{theorem}

The theorem improves the local limit law in \cite{kargin11b}, where it was
found that it holds for the window size $\eta^{\ast}\sim ( \log
N ) ^{-1/2}$. The optimal result is probably $\eta^{\ast}\sim
N^{-1+\varepsilon}$ with arbitrarily small positive $\varepsilon$, similar
to the case of classical Gaussian ensembles and the case of Wigner/sample
covariance matrices. The proof of Theorem \ref{Thmlocallaw} will be given
in Section~\ref{sectionlocallaw}.

\subsection{Largest eigenvalues of finite rank deformations of unitarily-invariant matrices}\label{sec14}

The largest eigenvalues of finite-rank deformations of Wigner matrices have
been recently received much attention and studied in \cite{peche06,maida07,feralpeche07,capitainedonati-martinferal09,capitainedonati-martinferalfevrier11,male11,pizzorenfrewsoshnikov11,renfrewsoshnikov12,knowlesyin12a,knowlesyin12b,capitaine12} and \cite
{peng12}. This study is closely connected to the study of spiked population
models in \cite{baikbenarouspeche05,baiksilverstein06,baiyao07,bloemendalvirag11} and \cite{bloemendalvirag11a}.

The idea that the subordination identities are useful in the context of
matrix deformations has first appeared in the work of Capitaine,
Donati-Martin, Feral and Fevrier (see \cite{capitainedonati-martinferalfevrier11} and \cite{capitaine12}). We will
use this idea to give a different proof for a result of Benaych-Georges and
Nadakuditi in \cite{BenaychGeorgesNadakuditi11}. They considered the
largest eigenvalue of $H_{N}=A_{N}+U_{N}B_{N}U_{N}^{\ast}$, where $A_{N}$
is a finite rank Hermitian matrix, and found a formula for the limits
of the
largest eigenvalues. (See also \cite{BenaychGeorgesGuionnetMaida11,BenaychGeorgesGuionnetMaida11a} and \cite{BenaychGeorgesNadakuditi12}
for further developments.)

Theorem \ref{PropRAestimate} allows us to obtain a different proof of
Benaych-Georges and Nadakuditi's result. While their method is based on
analysing the zeros of determinants of certain matrix-valued functions, our
method uses the singularities of the resolvent traces. In particular,
we use
the description that Theorem \ref{PropRAestimate} gives for the resolvent
behavior in the upper half-plane.

We consider the simplest case when matrix $A_{N}$ has rank one. The
ideas of
the proof can be applied similarly in the case when $A_{N}$ is a finite-rank
matrix with the rank fixed and $N$ approaching infinity.

Let $\rho_{\mu} ( \theta ) $ be the largest real solution
of the
equation $\theta m_{\mu} ( x ) +1=0$, and let $\lambda
_{1} (
X ) $ denote the largest eigenvalue of Hermitian matrix $X$.

\begin{theorem}
\label{theoremspikes}Let $H_{N}=A_{N}+U_{N}B_{N}U_{N}^{\ast}$ where $A_{N}$
is a rank-one Hermitian matrix with the eigenvalue $\theta_{0}>0$, and
$%
B_{N}$ is a Hermitian matrix with the empirical eigenvalue distribution
$\mu
_{B_{N}}$. Let $\lambda_{1} ( B_{N} ) \rightarrow L$ in
probability. Assume that matrices $B_{N}$ are uniformly bounded almost
surely and that $\mu_{B_{N}}$ weakly converges to a probability
measure $%
\mu$. Then
\[
\lambda_{1}^{ ( H_{N} ) }\rightarrow \cases{ \displaystyle
\rho_{\mu} ( \theta_{0} ), &\quad if $\rho _{\mu}
(\theta_{0} ) >L$,
\vspace*{3pt}\cr
L, &\quad otherwise,}
\]
where convergence is in probability.
\end{theorem}

The proof of Theorem \ref{theoremspikes} is based on the subordination-like
formula, which we will prove in Proposition \ref{propositionformulaspikes}:
%
\begin{equation}
\qquad \mathbb{E}m_{H_{N}} ( z ) =m_{B_{N}} ( z ) +
\frac{1}{N}%
\frac{m_{B_{N}}^{\prime} ( z ) }{m_{B_{N}} ( z
) } \biggl( \frac{1}{\theta m_{B_{N}} ( z ) +1}-1
\biggr) +O_{\eta
} \biggl( \frac{%
1}{N^{2}} \biggr),
\end{equation}
where $O_{\eta} ( N^{-2} ) $ denotes a function $f (
z ) $
such that $N^{2}\llvert  f ( z ) \rrvert \leq C
( \operatorname{Im}z ) ^{-k}$ for some $k>0$ and $C>0$. This formula explicitly shows
the correction term to the Stieltjes transform $m_{B_{N}} (
z ) $
that results from adding matrix $A_{N}$. In particular, this correction term
has an additional pole to the right of $L$ if and only if a zero of
$\theta
m_{B_{N}} ( z ) +1$ is located to the right of $L$. This implies
(after some additional work) that $\rho_{\mu} ( \theta
_{0} ) $
is the only possible limit point for the largest eigenvalue of $H_{N}$.

After the preprint of this paper has appeared, the method of subordination
functions was used in \cite{bbcf12} to generalize the results \cite{BenaychGeorgesNadakuditi11}. The main innovation in \cite{bbcf12} is that
the matrix $A_{N}$ is no longer required to be finite rank. It is only
required that it has sufficiently large fixed eigenvalues
(``spikes'').

Two standard examples in the deformation theory are Gaussian--Hermitian
matrices and Gaussian--Wishart matrices as $B_{N}$. In these examples,
Theorem~\ref{theoremspikes} gives the results in agreement with available in the
literature \cite{baikbenarouspeche05} and
\cite{capitainedonati-martinferal09}, with $\rho_{H} ( \theta
 )
=\theta+\sigma^{2}/\theta$ and $\rho_{W} ( \theta )
=\theta
+\lambda\theta/ ( \theta-1 ) $.

Another example, which seems to be new, is provided by random projection
matrices.

\begin{figure}

\includegraphics{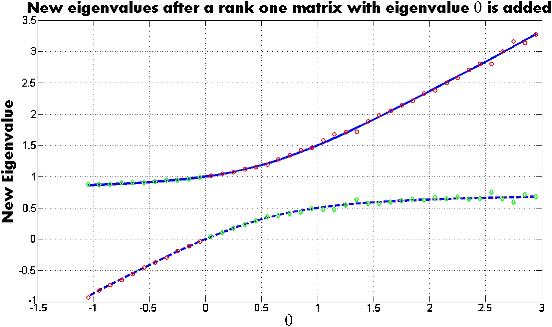}

\caption{New eigenvalues of the rank one perturbation of the projection
matrix model with $p=1/4$. Circles are eigenvalues outside of $[0,1]$,
diamonds are eigenvalues inside $[0,1]$. The solid line is $\protect\rho
(\theta)$, the dashed line is its conjugate.}\label{figurebernoulli}
\end{figure}

\begin{ex*}
Consider matrices $B_{N}=U_{N}P_{N}U_{N}^{\ast}$
where $%
P_{N}$ is a projection matrix of rank $p_{N}$. If $p_{N}/N\rightarrow p>0$
as $N\rightarrow\infty$, then the empirical eigenvalue distribution
of $%
B_{N}$ converges to the Bernoulli distribution $\mu_{b}=p\delta
_{1}+q\delta_{0}$, where $q=1-p$. One computes that
\[
\rho ( \theta ) =\frac{\theta}{2}+\frac{1}{2} \bigl( 1+
\sqrt{%
 ( 1+\theta ) ^{2}-4q\theta} \bigr),
\]
and this is the limit of the largest eigenvalue for the matrices $%
A_{N}+B_{N} $ when $N\rightarrow\infty$. This formula for the limit
of the
largest eigenvalue is valid for all $\theta_{0}>0$.

In the context of this example, an interesting phenomenon is uncovered by
numerical evidence, which is not explained by Theorem \ref{theoremspikes}.
Namely, adding a rank one projection $A_{N}$ with eigenvalue $\theta$
results in a creation of \textit{two} ``new'' eigenvalues. (See Figure~\ref{figurebernoulli}
for a
numerical example with the size of the matrices fixed at $N=100$.) One new
eigenvalue is given by $\rho ( \theta ) $, and another one
by the
other solution of the equation $\theta m ( x ) +1=0$:
\[
\overline{\rho} ( \theta ) =\frac{\theta}{2}+\frac
{1}{2} \bigl(
1-%
\sqrt{ ( 1+\theta ) ^{2}-4q\theta} \bigr).
\]
\end{ex*}

\subsection{Brief overview}\label{sec15}

In this paper, we consider the resolvent of the matrix $H=A+B$, where
$A:=V%
\widetilde{A}V^{\ast}$ and $B:=U\widetilde{B}U^{\ast}$, $\widetilde{A}$
and $\widetilde{B}$ are two $N$-by-$N$ Hermitian diagonal matrices,
and $U$
and $V$ are two $N$-by-$N$ random independent uniformly distributed unitary
matrices.

We have showed that there exist two functions $\omega_{A} (
z ) $
and $\omega_{B} ( z ) $ that depend only on the sets of
eigenvalues of $A$ and $B$ and have the following properties:
\begin{longlist}[(iii)]
\item[(i)] $\omega_{A} ( z ) $ and\vspace*{2pt} $\omega_{B} ( z ) $ are
analytic in $\mathbb{C}^{+}$;

\item[(ii)] if $\operatorname{Im}z\gg N^{-1/5}$, then
\[
\min \bigl\{ \operatorname{Im}\omega_{A} ( z ),\operatorname{Im}\omega
_{B} ( z ) \bigr\} \geq\operatorname{Im}z-\frac{c}{N (
\operatorname{Im}%
z ) ^{7}}.
\]
Moreover,\vspace*{2pt} if $\mu_{A}=\mu_{B}$, then $\min \{ \operatorname{Im}\omega
_{A} ( z ),\operatorname{Im}\omega_{B} ( z )
\} \geq
\operatorname{Im}z$ for all $z\in\mathbb{C}^{+}$;

\item[(iii)] if $\operatorname{Im}z\gg N^{-1/7}$, then
\[
\mathbb{E}_{U}G_{H} ( z ) -G_{A} \bigl(
\omega_{B} ( z ) \bigr) =O \biggl( \frac{1}{N\eta^{6}} \biggr)
\]
and
\[
\mathbb{E}m_{H} ( z ) -m_{A} \bigl( \omega_{B}
( z ) \bigr) =O \biggl( \frac{1}{N^{2}\eta^{6}} \biggr),
\]
and similar estimates hold for $\mathbb{E}_{V}G_{H} ( z )
-G_{B} ( \omega_{A} ( z )  ) $ and $\mathbb
{E}m_{H} (
z ) -m_{B} ( \omega_{A} ( z )  ) $.

This can be thought of as a subordination property for the resolvent of the
sum $A+B$ with respect to resolvents of $A$ and $B$.

We have used the subordination property to show that the localization length
of eigenvectors is greater than $N^{\kappa}$, where $\kappa=1/7$. The
probable actual localization length is $O ( N ) $.

Next, we have showed that a local limit law holds for the empirical
eigenvalue measure $\mu_{H_{N}}$ with the window length $N^{-1/7}$. This
result improves over the result in \cite{kargin11b}. However, it is still
far from the probable optimal result with the window length $%
N^{-1+\varepsilon}$.

Finally, by using our results about subordination we studied the rank-one
deformations of unitarily-invariant random matrices, and derived explicit
formulas for the limit of their largest eigenvalues.

The rest of the paper is organized as follows. Section \ref{sectionsubordination} is devoted to the proof of Theorems \ref{Propsubordination} and \ref{PropRAestimate} regarding the
subordination. Section~\ref{sectiondelocalization} is about delocalization
of eigenvectors (Theorem \ref{Thmdelocalization}). Section \ref{sectionlocallaw} proves Theorem \ref{Thmlocallaw} about the local law
for eigenvalues. Section~\ref{sectionspikes} proves Theorem \ref{theoremspikes} about rank-one deformations of unitarily-invariant
ensembles. And three appendices contain various auxiliar results.
\end{longlist}

\section{Approximate subordination}\label{sec2}
\label{sectionsubordination}

Before we start the proof of Theorems \ref{Propsubordination} and
\ref{PropRAestimate}, note that the definitions imply the following useful
identity:
%
\begin{equation}
\omega_{A} ( z ) +\omega_{B} ( z ) =z-\frac
{1}{\mathbb{E}%
m_{H} ( z ) }.
\label{identitysumomegas}
\end{equation}
Indeed,
\[
\omega_{A} ( z ) +\omega_{B} ( z ) =2z-\frac
{\mathbb{E}%
 [ f_{A} ( z ) +f_{B} ( z )  ]
}{\mathbb{E}%
m_{H} ( z ) }
\]
and
\begin{eqnarray*}
f_{A} ( z ) +f_{B} ( z ) &=&N^{-1}\operatorname{Tr}
\bigl( ( A+B ) ( A+B-zI ) ^{-1} \bigr)
\\
&=&1+zN^{-1}\operatorname{Tr} \bigl( ( A+B-zI ) ^{-1} \bigr)
\\
&=&1+zm_{H} ( z ),
\end{eqnarray*}
which implies (\ref{identitysumomegas}).

Now we start proving Theorem \ref{Propsubordination}. First, write
%
\begin{equation}
\mathbb{E}_{U}G_{H} ( z ) =G_{A} \bigl(
\omega_{B} ( z ) \bigr) +R_{A} ( z ). \label{formulasubordination}
\end{equation}
The error term in subordination formula (\ref{formulasubordination})
can be
written as follows:
%
\begin{equation}
R_{A} ( z ):=\frac{1}{\mathbb{E}m_{H}}(A-zI)G_{A} \bigl( \omega
_{B} ( z ) \bigr) \mathbb{E}_{U}\Delta_{A},
\label{13154845497151451}
\end{equation}
where
\[
\Delta_{A}:=- ( m_{H}-\mathbb{E}m_{H} )
G_{H}-G_{A} ( f_{B}-%
\mathbb{E}f_{B} ) G_{H}.
\]
In order to derive formulae (\ref{formulasubordination}) and (\ref{13154845497151451}), one starts by calculating $dG_{t}/dt$ where $%
G_{t}= ( A+e^{iXt}Be^{-iXt} ) ^{-1}$ and $X$ is an Hermitian
matrix. Since $B$ has a rotationally invariant distribution, hence
$\mathbb{E%
}_{U} ( dG_{t}/dt ) =0$, and one can find by using different
generator matrices $X$ that this implies that $\mathbb{E}_{U} (
G_{H}\otimes BG_{H} ) = \mathbb{E}_{U} ( G_{H}B\otimes
G_{H} ) $. After taking the trace over the first component of the
tensor product, one gets $\mathbb{E}_{U} ( m_{H}BG_{H} ) =
\mathbb{E}_{U} ( f_{B}G_{H} ) $. This can be rewritten as
$\mathbb{%
E}_{U} ( m_{H}G_{H} ) =G_{A}\mathbb{E}_{U} (
m_{H}I-f_{B}G_{H} ) $. Next, one writes $\mathbb{E}_{U} (
m_{H}G_{H} ) =\mathbb{E} ( m_{H} ) \mathbb
{E}_{U} (
G_{H} ) +e_{1}$ and $\mathbb{E}_{U} ( f_{B}G_{H} )
=\mathbb{E}%
 ( f_{B} ) \mathbb{E}_{U} ( G_{H} ) +e_{2}$,
where $e_{1}$
and $e_{2}$ are error terms. After substituting these expressions, one can
manipulate the previous identity so that $\mathbb{E}_{U} (
G_{H} ) $
is on the left-hand side and everything else is on the right-hand side. The
resulting expression is equivalent to (\ref{formulasubordination})
with the
error term given by (\ref{13154845497151451}). See Appendix \ref{secB} for a more
complete derivation, and \cite{pasturvasilchuk00} or proof of Theorem
7 in~\cite{kargin11b} for details.

We can also rewrite formula (\ref{formulasubordination}) as follows:
\[
\mathbb{E}_{U}G_{H}=G_{A} \bigl(
\omega_{B} ( z ) \bigr) \biggl( I+%
\frac{1}{\mathbb{E}m_{H}} ( A-zI
) \mathbb{E}_{U}\Delta _{A} \biggr).
\]
Hence,
\[
( \mathbb{E}_{U}G_{H} ) ^{-1}= \biggl( I+
\frac{1}{\mathbb
{E}m_{H}}%
 ( A-zI ) \mathbb{E}_{U}
\Delta_{A} \biggr) ^{-1} \bigl( A-\omega _{B} ( z )
I \bigr)
\]
and
%
\begin{eqnarray}\label{formulatrueinversion}
\omega_{B} ( z ) I &=&- ( \mathbb{E}_{U}G_{H} )
^{-1}+A
\nonumber\\[-8pt]\\[-8pt]
&&{} + \biggl[ \biggl( I+\frac{1}{\mathbb{E}m_{H}} ( A-zI )
\mathbb{E}%
_{U}\Delta_{A} \biggr) ^{-1}-I
\biggr] \bigl( A-\omega_{B} ( z ) I \bigr). \nonumber
\end{eqnarray}
Let us consider the first two terms in this expression. Later, we are going
to show that the third term is small.

Define
\[
\Omega_{B} ( z,A ):=- \bigl( \mathbb{E}_{U}G_{H}
( z ) \bigr) ^{-1}+A.
\]
The matrix function $\Omega_{B} ( z,A ) $ has a property
which is
similar to the subordination property.

\begin{lemma}
\label{propositionsubordination1} Let $\lambda ( z ) $ be an
eigenvalue of $\Omega_{B} ( z,A ) $. Then $\operatorname{Im}\lambda
 ( z ) \geq\operatorname{Im}z$.
\end{lemma}

\begin{pf}
Let $z=x+i\eta$ and $\eta>0$. Then every matrix
$ (
H-x-i\eta ) ^{-1}$ is normal and its eigenvalues are on the
border of
a disc $D_{\eta}$ with the center at $i/ ( 2\eta ) $ and the
radius equal to $1/ ( 2\eta ) $. Hence, by Lemma \ref{lemmaconvcombeigs} in Appendix~\ref{secC}, the eigenvalues of $\mathbb{E}%
_{U} ( H-x-i\eta ) ^{-1}$ belong to the disc $D_{\eta}$. It
follows that eigenvalues of $- [ \mathbb{E}_{U} ( H-x-i\eta
 )
^{-1} ] ^{-1}$ are in $\mathbb{H}_{\eta}= \{ w\dvtx \operatorname{Im}w\geq
\eta \} $.

If we take the basis in which $A$ is diagonal, then $ ( \mathbb{E}
_{U}G_{H} ( z )  ) ^{-1}= [ \mathbb{E}_{U}
( H-x-i\eta
 ) ^{-1} ] ^{-1}$ is diagonal by Lemma \ref
{PropEGdiagonal} in
Appendix~\ref{secC}. Since $A$ is Hermitian, therefore its eigenvalues are real.
Hence, the imaginary parts of eigenvalues of $\Omega_{B} (
z,A ) $
coincide with imaginary parts of eigenvalues of $- [ \mathbb{E}%
_{U} ( H-x-i\eta ) ^{-1} ] ^{-1}$, and we arrive at
the claim
of the lemma.
\end{pf}

Now we are going to estimate the size of the third term in the right-hand
side in~(\ref{formulatrueinversion}). First, we estimate the size of
$%
\mathbb{E}_{U}\Delta_{A}$. We use concentration inequalities.

\begin{lemma}
\label{lemmaDeltaAestimate} Assume that $\eta:=\operatorname{Im}z\in
 (
0,1 ) $ and $\llvert \operatorname{Re}z\rrvert \leq K (
A,B )
$. Then
\[
\mathbb{E}_{U}\Delta_{A} ( z ) =O \biggl(
\frac{1}{\eta
^{4}N}%
 \biggr).
\]
\end{lemma}

\begin{pf}
Since$\llVert  G_{H}\rrVert \leq1/\eta$,
hence by
using Lemma \ref{lemmadiffGiiEGii} in Appendix \ref{secD}, we obtain
\begin{eqnarray*}
\mathbb{P} \bigl\{ \bigl\llVert \bigl( m_{H} ( z ) -\mathbb
{E}%
m_{H} ( z ) \bigr) G_{H}\bigr\rrVert \geq
\delta/\eta \bigr\} &\leq&\exp \biggl[ -c\frac{\delta^{2}\eta^{4}}{\llVert  B\rrVert ^{2}%
}N^{2} \biggr],
\\
\mathbb{P} \bigl\{ \bigl\llVert G_{A} \bigl( f_{B} ( z ) -
\mathbb{E}%
f_{B} ( z ) \bigr) G_{H}\bigr\rrVert
\geq\delta/\eta ^{2} \bigr\} &\leq&\exp \biggl[ -c\frac{\delta^{2}\eta^{4}}{\llVert  B\rrVert ^{2}%
}N^{2}
\biggr].
\end{eqnarray*}
Set $\varepsilon=\delta/\eta$ and $\varepsilon=\delta/\eta^{2}$
in the
first and the second inequalities, respectively, and use the triangle
inequality for norms in order to obtain that
\begin{eqnarray*}
\mathbb{P} \bigl\{ \bigl\llVert \Delta_{A} ( z ) \bigr\rrVert \geq
\varepsilon \bigr\} &\leq&\exp \biggl[ -\frac{c\varepsilon
^{2}N^{2}}{%
\llVert  B\rrVert ^{2}}\min \bigl\{
\eta^{6},\eta^{8} \bigr\} \biggr]
\\
&\leq&\exp \bigl[ -c\varepsilon^{2}\eta^{8}N^{2}
\bigr].
\end{eqnarray*}
Next, note that $\llVert \mathbb{E}_{U}\Delta_{A}\rrVert
\leq
\mathbb{E}_{U}\llVert \Delta_{A}\rrVert $ by the convexity
of norm,
and $\mathbb{E}_{U}\llVert \Delta_{A}\rrVert $ can be
estimated by
using the equality $\mathbb{E}X=\int_{0}^{\infty} ( 1-\mathcal
{F}%
_{X} ( t )  ) \,dt$, valid for every positive random
variable $X$
and its cumulative distribution function $\mathcal{F}_{X} (
t ) $.
In our case, we obtain
\[
\mathbb{E}_{U}\llVert \Delta_{A}\rrVert \leq\int
_{0}^{\infty}\exp%
 \bigl[ -ct^{2}
\eta^{8}N^{2} \bigr] \,dt=\frac{c^{\prime}}{N\eta^{4}}.
\]\upqed
\end{pf}

Next, Lemma \ref{lemmaimmzlowerbound} in Appendix \ref{secD} says that
$ (
\mathbb{E}m_{H} ( z )  ) ^{-1}\leq c/\eta$. Hence,
Lemma~\ref{lemmaDeltaAestimate} implies that
%
\begin{equation}
\biggl\llVert \frac{1}{\mathbb{E}m_{H}} ( A-z ) \mathbb {E}\Delta _{A}\biggr
\rrVert \leq\frac{c}{\eta^{5}N}, \label{psiestimate}
\end{equation}
where $c>0$ depends only on $K$ and $R$.

It is easy to prove that if $\llVert  X\rrVert \leq\varepsilon<1/2$,
then $\llVert  ( I+X ) ^{-1}-I\rrVert \leq
2\varepsilon$.
In particular, for all $N\gg\eta^{-5}$, we have
\[
\biggl\llVert \biggl( I+\frac{1}{\mathbb{E}m_{H}} ( A-z ) \mathbb{E}%
\Delta_{A} \biggr) ^{-1}-I\biggr\rrVert \leq
\frac{c}{\eta^{5}N}.
\]

Next, note that by definition $\omega_{B} ( z ) =z-\mathbb
{E}%
f_{B} ( z ) /\mathbb{E}m_{H} ( z ) $. From Lemma
\ref{lemmaimmzlowerbound}, $\llvert  ( \mathbb{E}m_{H} (
z )
 ) ^{-1}\rrvert <c/\eta$. In addition,
\[
\bigl\llvert \mathbb{E}f_{B} ( z ) \bigr\rrvert =\biggl\llvert
\mathbb{E}%
\frac{1}{N}\operatorname{Tr} \biggl( B\frac{1}{H-z}
\biggr) \biggr\rrvert \leq \llVert B\rrVert \mathbb{E} \biggl( \biggl\llVert
\frac
{1}{H-z}\biggr\rrVert \biggr) \leq c\frac{1}{\eta}.
\]
Hence, $\llvert \omega_{B} ( z ) -z\rrvert \leq
c/\eta^{2}$.
It follows that
%
\begin{equation}
\biggl\llVert \biggl( \biggl( I+\frac{1}{\mathbb{E}m_{H}} ( A-z ) \mathbb{E}
\Delta_{A} \biggr) ^{-1}-I \biggr) \bigl( A-\omega
_{B} ( z ) \bigr) \biggr\rrVert \leq\frac{c}{N\eta^{7}}. \label{errortermestimate}
\end{equation}

\begin{lemma}
\label{lemmaperturbationnormaloper}Let $\Omega$ be a diagonal
matrix and
$R$ be an arbitrary matrix. Then for every eigenvalue $\widehat
{\lambda}%
_{i} $ of $\Omega+R$, there exists an eigenvalue $\lambda_{i}$ of
$\Omega$
such that $\llvert \widehat{\lambda}_{i}-\lambda_{i}\rrvert
\leq
\llVert  R\rrVert $.
\end{lemma}

(See Theorem 6.3.2 on page 365 in \cite{hornjohnson85}.)

Formulae (\ref{formulatrueinversion}) and (\ref
{errortermestimate}), and
Lemmas \ref{propositionsubordination1} and \ref{lemmaperturbationnormaloper} imply that
\[
\operatorname{Im} \bigl( \omega_{B} ( z ) \bigr) \geq\operatorname{Im}z-
\frac{%
c}{N ( \operatorname{Im}z ) ^{7}}.
\]
This completes the proof of Theorem \ref{Propsubordination}.

\begin{pf*}{Proof of Theorem \ref{PropRAestimate}}
Since $N\gg1/\eta^{7}$,
Theorem \ref{Propsubordination} implies that $\operatorname{Im} (
\omega
_{B} ( z )  ) \geq\eta/2$. Hence, $\llVert
G_{A} (
\omega_{B} ( z )  ) \rrVert \leq c/\eta$. Since
\[
R_{A} ( z ) =G_{A} \bigl( \omega_{B} ( z )
\bigr) \frac{1%
}{\mathbb{E}m_{H}} ( A-z ) \mathbb{E}_{U}\Delta_{A},
\]
we can use (\ref{psiestimate}) in order to obtain
\[
\bigl\llVert R_{A} ( z ) \bigr\rrVert \leq\bigl\llVert
G_{A} \bigl( \omega_{B} ( z ) \bigr) \bigr\rrVert \biggl
\llVert \frac
{1}{\mathbb{E}%
m_{H}} ( A-z ) \mathbb{E}_{U}\Delta_{A}
\biggr\rrVert \leq \frac{c}{%
N\eta^{6}},
\]
which yields the first point of Theorem \ref{PropRAestimate}.

In order to estimate the error term in the second part of the theorem, we
note that by definition of $R_{A}$ it is enough to show that
\[
\mathbb{E}\frac{1}{N}\operatorname{Tr} \bigl[ (A-z)G_{A} \bigl(
\omega _{B} ( z ) \bigr) \Delta_{A} \bigr] =O \biggl(
\frac{1}{\eta
^{6}N^{2}} \biggr).
\]

By using the definition of $\Delta_{A}$, we can write the modulus of the
expression on the left-hand side as follows:
\[
\bigl\llvert \mathbb{E} ( m_{H}-\mathbb{E}m_{H} ) (
\varphi-%
\mathbb{E}\varphi ) +\mathbb{E} ( f_{B}-\mathbb
{E}f_{B} ) ( \psi-\mathbb{E}\psi ) \bigr\rrvert,
\]
where
\[
\varphi:=\frac{1}{N}\operatorname{Tr} \bigl[ ( A-z ) G_{A} \bigl(
\omega_{B} ( z ) \bigr) G_{H} ( z ) \bigr]
\]
and
\[
\psi:=\frac{1}{N}\operatorname{Tr} \bigl[ ( A-z ) G_{A} \bigl(
\omega _{B} ( z ) \bigr) G_{A} ( z ) G_{H} ( z )
\bigr].
\]
By the Cauchy--Schwarz inequality, we estimate this from above by
%
\begin{equation}
\sqrt{\operatorname{Var}(m_{H})\operatorname{Var} ( \varphi ) }+\sqrt{
\operatorname{Var}(f_{B})\operatorname{Var} ( \psi ) }. \label{errorforStieltjes}
\end{equation}
By applying Lemma \ref{lemmadiffGiiEGii} from Appendix \ref{secD} and the estimate
$\llVert  G_{A} ( \omega_{B} ) \rrVert \leq c/\eta
$ in
order to bound the variances, we find\vspace*{1pt} that for sufficiently large $N$, $
\operatorname{Var} ( m_{H} ) =O ( \eta^{-4}N^{-2} )
$, $\operatorname{Var} ( \varphi ) =O ( \eta^{-6}N^{-2} ) $,
$\operatorname{Var}%
 ( f_{B} ) =O ( \eta^{-4}N^{-2} ) $ and
$\operatorname{Var}%
 ( \psi ) =O ( \eta^{-8}N^{-2} ) $. Hence, the
expression in (\ref{errorforStieltjes}) is smaller than $c/ (
\eta
^{6}N^{2} ) $, provided that $z$ is in the region where $\omega
_{B} ( z ) $ increases the imaginary part. This completes
the proof
of the second part of the theorem.

The\vspace*{1pt} third part of Theorem \ref{PropRAestimate} immediately follows from
the first part and Lemma \ref{lemmadiffGiiEGii} in Appendix \ref{secD} if we
take $%
\eta\gg N^{-1/7}$.

Indeed, if $\operatorname{Im}z=\eta\gg N^{-1/7}\gg N^{-1/6}$, then the
first part
of Theorem \ref{PropRAestimate} implies that $\llvert  (
\mathbb{E}%
_{U}G_{H} ( z )  ) _{ij}- ( G_{A} ( \omega
_{B} (
z )  )  ) _{ij}\rrvert \leq\delta/2$ for all
sufficiently large $N$. For these $N$, we have
\begin{eqnarray*}
&& \mathbb{P}_{U} \bigl\{ \bigl\llvert \bigl( G_{H} ( z )
\bigr) _{ij}- \bigl( G_{A} \bigl( \omega_{B} ( z )
\bigr) \bigr) _{ij}\bigr\rrvert \geq\delta \bigr\}
\\
&&\qquad \leq
\mathbb{P}_{U} \bigl\{ \bigl\llvert \bigl( G_{H} ( z )
\bigr) _{ij}- \bigl( \mathbb{E}%
_{U}G_{H}
( z ) \bigr) _{ij}\bigr\rrvert \geq\delta /2 \bigr\}
\\
&&\qquad \leq\exp \biggl( -\frac{c\delta^{2}\eta^{4}}{\llVert  B\rrVert ^{2}%
}N \biggr) \leq\exp \bigl( -c
\delta^{2}N^{3/7} \bigr).
\end{eqnarray*}\upqed
\end{pf*}

\section{Delocalization}\label{sec3}
\label{sectiondelocalization}

The essential part of the proof is to show that $\omega_{A_{N}}$ and $%
\omega_{B_{N}}$ are close to $\omega_{\alpha}$ and $\omega_{\beta}$,
respectively. Namely, let
%
\begin{equation}
r ( z ):=\max\bigl\{\bigl\llvert r_{A} ( z ) \bigr\rrvert,\bigl
\llvert r_{B} ( z ) \bigr\rrvert \bigr\} \label{definitionr}
\end{equation}
and
%
\begin{equation}
s ( A,B ):=\max\bigl\{d_{L}(\mu_{A},
\mu_{\alpha}),d_{L}(\mu _{B},\mu_{\beta})
\bigr\}. \label{definitions}
\end{equation}

\begin{proposition}
\label{propdiffsubordination} Assume that a pair of probability
measures $%
 ( \mu_{\alpha},\mu_{\beta} ) $ is smooth in a closed
interval $%
I$. Then for some positive $\overline{r}$, $\overline{s}$ and
$\overline{%
\eta}$, if $r ( z ) \leq\overline{r}$, $s (
A,B ) \leq
\overline{s}$, $\operatorname{Re}z\in I$ and $\operatorname{Im}z\in(0,\overline
{\eta}]$, then
\[
\max \bigl( \bigl\llvert \omega_{\alpha} ( z ) -\omega _{A} ( z
) \bigr\rrvert,\bigl\llvert \omega_{\beta} ( z ) -\omega _{B}
( z ) \bigr\rrvert \bigr) =O ( r+s ),
\]
where the constant in the $O$-term may depend on the pair $ ( \mu
_{\alpha},\mu_{\beta} ) $ and on $\max \{ \llVert
A\rrVert,\llVert  B\rrVert  \} $.
\end{proposition}

Let us postpone the proof and show how this result implies Theorem \ref{Thmdelocalization}.

\begin{pf*}{Proof of Theorem \ref{Thmdelocalization}}
Let $N$ be the
size of
matrices $A$ and $B$ and assume that $N$ is sufficiently large so that
$\max
 \{ d_{L}(\mu_{A},\mu_{\alpha}),d_{L}(\mu_{B},\mu_{\beta
}) \}
\leq s<\overline{s}$. By definition, $r_{A} ( z ) =\mathbb
{E}%
m_{H} ( z ) -m_{A} ( \omega_{B} ( z )
 ) $,
hence the second part of Theorem \ref{PropRAestimate} says that if
$N\gg
\eta^{-7}$, then $\llvert  r_{A} ( z ) \rrvert
=O (
\frac{1}{N^{2}\eta^{6}} ) $. A similar bound holds for $\llvert
r_{B} ( z ) \rrvert $. Hence, we can take $r=O (
\frac{1}{%
N^{2}\eta^{6}} ) $ in Proposition \ref{propdiffsubordination} and
conclude that
%
\begin{equation}
\omega_{\beta} ( x+i\eta ) -\omega_{B} ( x+i\eta ) =O \biggl(
\frac{1}{N^{2}\eta^{6}}+s \biggr). \label{omegaBandomegabeta}
\end{equation}

Hence, if $s$ is sufficiently small, then $\operatorname{Im}\omega_{B} (
z ) \geq c>0$ for all $z$ with $\operatorname{Re}z\in I$ and
$cN^{-2/6}\leq
\operatorname{Im}z\leq\overline{\eta}$. It follows that $ [
G_{A} (
\omega_{B} ( z )  )  ] _{kk}$ is bounded, say, $
\llvert  [ G_{A} ( \omega_{B} ( z )  )
 ]
_{kk}\rrvert <C$. By using the third part of Theorem \ref{PropRAestimate}, we find that
\[
\mathbb{P} \bigl\{ \bigl\llvert \bigl[ G_{H} ( z ) \bigr]
_{kk}\bigr\rrvert \geq C+\delta \bigr\} \leq\exp \bigl( -c\delta
^{2}N^{3/7} \bigr).
\]
Now let $ \{ v_{a} \} _{a=1}^{N}$ denote an orthonormal
basis of
eigenvectors of $H$ and let $\lambda_{a}$ be the corresponding eigenvalues.
Let $v_{a} ( j ) $ denote the $j$th component of vector
$v_{a}$ in
the standard basis. Since
\[
G_{H} ( z ) =\sum_{a=1}^{N}
\frac{\llvert  v_{a}
\rangle
 \langle v_{a}\rrvert }{\lambda_{a}-z},
\]
hence
\[
\operatorname{Im}G_{kk} ( x+i\eta ) =\sum_{a=1}^{N}
\frac{\eta
\llvert
v_{a} ( k ) \rrvert ^{2}}{ ( \lambda_{a}-x )
^{2}+\eta^{2}}.
\]

Let us set $x=\lambda_{a}$ for a particular value of $a$, then
\[
\operatorname{Im}G_{kk} ( \lambda_{a}+i\eta ) \geq
\frac{\llvert
v_{a} ( k ) \rrvert ^{2}}{\eta},
\]
and, therefore,
%
\begin{equation}
\bigl\llvert v_{a} ( k ) \bigr\rrvert ^{2}\leq\eta\operatorname{Im}%
G_{kk} ( \lambda_{a}+i\eta ) \leq C\eta\leq
CN^{-1/7}\log N. \label{boundeigenvector}
\end{equation}\upqed
\end{pf*}

Before starting the proof of Proposition \ref
{propdiffsubordination}, let
us exclude $m ( z ) $ from the free probability system (\ref
{system}%
):
%
\begin{eqnarray}\label{systemfreelimit}
m_{\alpha} \bigl( \omega_{\beta} ( z ) \bigr) +\frac
{1}{\omega
_{\alpha} ( z ) +\omega_{\beta} ( z ) -z}
&=&0,
\nonumber\\[-8pt]\\[-8pt]
m_{\beta} \bigl( \omega_{\alpha} ( z ) \bigr) +\frac
{1}{\omega
_{\alpha} ( z ) +\omega_{\beta} ( z ) -z}
&=&0.\nonumber
\end{eqnarray}
A similar system can be written in the matrix case for $\omega
_{A} (
z ) $ and $\omega_{B} ( z ) $:
%
\begin{eqnarray}\label{systemmatrix}
m_{A} \bigl( \omega_{B} ( z ) \bigr) +\frac{1}{\omega
_{A} (
z ) +\omega_{B} ( z ) -z}
&=&-r_{A} ( z ),
\nonumber\\[-8pt]\\[-8pt]
m_{B} \bigl( \omega_{A} ( z ) \bigr) +\frac{1}{\omega
_{A} (
z ) +\omega_{B} ( z ) -z}
&=&-r_{B} ( z ),
\nonumber
\end{eqnarray}
where $r_{A} ( z ):=N^{-1}\operatorname{Tr} ( R_{A} (
z )
 ) =N^{-1}\operatorname{Tr} ( \mathbb{E}R_{A} ( z )
 ) $, $%
r_{B} ( z ):=N^{-1}\operatorname{Tr} ( R_{B} ( z
)  )
=N^{-1}\operatorname{Tr} ( \mathbb{E}R_{B} ( z )  )
$. Here, $%
R_{A} ( z ) $ is defined in (\ref{formulasubordination}),
and $%
R_{B} ( z ):=\mathbb{E}_{V}G_{H} ( z )
-G_{B} (
\omega_{A} ( z )  ) $.

The proof of Proposition \ref{propdiffsubordination} is done by an
application of the Kantorovich--Newton method that allows us to study
how the
perturbation of the system for $\omega_{\alpha}$ and $\omega_{\beta}$
affects the solution. The role of Theorem \ref{PropRAestimate} in the
proof is to ensure that the size of the perturbation is small.

Let us briefly recall the Newton--Kantorovich method of successive
approximations \cite{kantorovich48}. The method is quite general and works
for perturbabions of maps acting on Banach spaces. We will use it for the
maps defined on pairs of functions $w_{1} ( z ) $,
$w_{2} (
z ) $ which are holomorphic in a compact domain $\Omega$. However,
since the maps can be considered for every $z$ separately, we will\vspace*{2pt}
essentially consider them as maps from $\mathbb{C}^{2}$ to $\mathbb{C}^{2}$
with the norm $\llVert  ( w_{1},w_{2} ) \rrVert
= (
\llvert  w_{1}\rrvert ^{2}+\llvert  w_{2}\rrvert
^{2} )
^{1/2}$.

The general setup is as follows. Let $F ( w ) =0$ be a nonlinear
functional equation where $F$ is a nonlinear operator that sends elements
of a Banach space $W$ to itself. Let $F$ be twice differentiable, and
assume that in a neighborhood of a point $w_{0}$ the operator
$F^{\prime
} ( w ) $ has an inverse $ [ F^{\prime} ( w
)  ]
^{-1}\in L ( W ) $ where $L ( W ) $ denotes the
space of
bounded linear operators from $W$ to $W$. Consider the iterations
\[
w_{n+1}=w_{n}- \bigl[ F^{\prime} ( w_{n}
) \bigr] ^{-1}F ( w_{n} ).
\]

The Kantorovich theorem (i) gives the sufficient conditions for the
convergence of this process to a solution $w^{\ast}$ of equation
$F (
w ) =0$, (ii) estimates the speed of convergence, and (iii) estimates
the distance of the solution $w^{\ast}$ from the initial point
$w_{0}$. We
give the statement of the theorem omitting the claim about the speed of
convergence, which is not important for us.

\begin{theorem}[(Kantorovich)]\label{theoremKantorovich}Suppose that the following conditions hold:%
\begin{longlist}[(iii)]
\item[(i)] for an initial approximation $w_{0}$, the operator $F^{\prime
} (
w_{0} ) $ possesses an inverse operator $\Gamma_{0}= [
F^{\prime
} ( w_{0} )  ] ^{-1}$ whose norm has the following
estimate: $%
\llVert \Gamma_{0}\rrVert \leq C_{0}$,

\item[(ii)] $\llVert \Gamma_{0}F ( w_{0} ) \rrVert \leq
\delta
_{0}$,

\item[(iii)] the second derivative $F^{\prime\prime} ( w ) $ is bounded
in the domain determined by inequality (\ref{inequalityKantorovich}) below,
namely, $\llVert  F^{\prime\prime} ( w ) \rrVert
\leq M$,%

\item[(iv)] the constants $C_{0},\delta_{0},M$ satisfy the relation $%
h_{0}=C_{0}\delta_{0}M\leq1/2$.

Then equation $F ( w ) =0$ has a solution $w^{\ast}$, which lies
in a neighborhood of $w_{0}$ determined by the inequality
%
\begin{equation}
\llVert w-w_{0}\rrVert \leq\frac{1-\sqrt
{1-2h_{0}}}{h_{0}}\delta _{0},
\label{inequalityKantorovich}
\end{equation}
and the successive approximations $w_{n}$ of the Newton method converge
to $%
w^{\ast}$.
\end{longlist}
\end{theorem}

\begin{pf*}{Proof of Proposition \ref{propdiffsubordination}}
We want to
prove that the solutions of systems~(\ref{systemfreelimit}) and
(\ref{systemmatrix}) are close to each other in a certain region of
$\mathbb{C}%
^{+}$.

Write system (\ref{systemmatrix}) as $F ( w ) =0$, where
\[
F\dvtx \pmatrix{ w_{1}
\vspace*{3pt}\cr
w_{2}}\rightarrow \pmatrix{ (
w_{1}+w_{2}-z ) ^{-1}+m_{A} (
w_{2} )+r_{A} (z )
\cr
( w_{1}+w_{2}-z
) ^{-1}+m_{B} ( w_{1} )+r_{B} (z )}.
\]

Theorem \ref{theoremKantorovich} requires estimating three norms, $%
\llVert \Gamma_{0}\rrVert $, $\llVert \Gamma_{0}F (
w_{0} ) \rrVert $ and\break $\llVert  F^{\prime\prime} (
w ) \rrVert $. We start by estimating the norm $\llVert
F (
w_{0} ) \rrVert $ with $w_{0}= ( \omega_{\alpha} (
z ),\omega_{\beta} ( z )  ) $. We have
\[
\bigl\llVert F ( w_{0} ) \bigr\rrVert =\Biggl\llVert \lleft(
\begin{array} {c} m_{A} \bigl( \omega_{\beta} ( z )
\bigr) -m_{\alpha
} \bigl( \omega_{\beta} ( z ) \bigr)
+r_{A} ( z )
\\[3pt]
m_{B} \bigl( \omega_{\alpha} ( z ) \bigr) -m_{\beta
}
\bigl( \omega_{\alpha} ( z ) \bigr) +r_{B} ( z )%
\end{array} %
 \rright) \Biggr\rrVert.
\]

By assumption, $\llVert  ( r_{A} ( z ),r_{B}
( z )
 ) \rrVert \leq r$. To complete the estimate, we also need a lemma.

\begin{lemma}
\label{lemmaclosenessSttransforms} Let $m_{1} ( z ) $
and $%
m_{2} ( z ) $ denote the Stieltjes transforms of measures
$\mu_{1}$
and $\mu_{2}$, respectively. Let $d_{L} ( \mu_{1},\mu_{2}
) =s$
and $z=x+i\eta$, where $\eta>0$. Then:
\begin{longlist}[(a)]
\item[(a)] $\llvert  m_{1} ( z ) -m_{2} ( z ) \rrvert
<cs\eta^{-1} \max \{ 1,\eta^{-1} \} $ where $c>0$ is a numeric
constant, and

\item[(b)] $\llvert \frac{d^{r}}{dz^{r}} ( m_{1} ( z )
-m_{2} (
z )  ) \rrvert <c_{r}s\eta^{-1-r}\max \{ 1,\eta
^{-1} \} $ where $c>0$ are numeric constants.
\end{longlist}
\end{lemma}

This lemma was proved as Lemma 2.2 in \cite{kargin13a}.

By assumption of smoothness on the closed inteval $I$, we know that
$\operatorname{Im}\omega_{\alpha} ( x ) $ and $\operatorname{Im}\omega_{\beta
} (
x ) $ are uniformly bounded away from zero on $I$. Moreover,
since $%
\omega_{\alpha}(z) $ and $\omega_{\beta}(z)$ are continuous in a
rectangle $R_{\varepsilon}:= \{ z=x+i\eta|x\in I,0\leq\eta\leq
\varepsilon \} $ (again by assumption of smoothness in the
interval $I$%
), hence $\operatorname{Im}\omega_{\alpha} ( z ) $ and $\operatorname{Im}%
\omega_{\beta} ( z ) $ are uniformly bounded away from
zero on
this rectangle provided that $\varepsilon$ is sufficiently small. We will
use this fact repeatedly below.

In particular, together with Lemma \ref{lemmaclosenessSttransforms} this
implies that
\[
\lleft\llVert \lleft( %
\begin{array} {c} m_{A} \bigl(
\omega_{\beta} ( z ) \bigr) -m_{\alpha
} \bigl(
\omega_{\beta} ( z ) \bigr)
\\[3pt]
m_{B} \bigl( \omega_{\alpha} ( z ) \bigr) -m_{\beta
}
\bigl( \omega_{\alpha} ( z ) \bigr)%
\end{array} %
 \rright) \rright\rrVert \leq cs
\]
on $R_{\varepsilon}$.

Hence, $\llVert  F ( w_{0} ) \rrVert $ is bounded by
$O (
r+s ) $ uniformly for every point $z\in R_{\varepsilon}$.

The next step is to estimate the norm of the inverse derivative. We compute
\[
F^{\prime}=\lleft( %
\begin{array} {cc} - (
w_{1}+w_{2}-z ) ^{-2} & - ( w_{1}+w_{2}-z
) ^{-2}+m_{A}^{\prime} ( w_{2} )
\\[3pt]
- ( w_{1}+w_{2}-z ) ^{-2}+m_{B}^{\prime}
( w_{1} ) & - ( w_{1}+w_{2}-z )
^{-2}%
\end{array} %
 \rright).
\]
The determinant of this matrix is
\[
\bigl[ m_{A}^{\prime} ( w_{2} ) +m_{B}^{\prime}
( w_{1} ) \bigr] ( w_{1}+w_{2}-z )
^{-2}-m_{A}^{\prime} ( w_{2} )
m_{B}^{\prime} ( w_{1} ).
\]

Recall that condition (B) in the assumption of smoothness requires that
%
\begin{equation}
k_{\mu}(x):=\frac{1}{m_{\mu_{\alpha}}^{\prime}(\omega_{\beta
}(x))}+%
\frac{1}{m_{\mu_{\beta}}^{\prime}(\omega_{\alpha}(x))}-\bigl(
\omega _{\alpha
}(x)+\omega_{\beta}(x)-x\bigr)^{2}\neq0.
\end{equation}
By continuity of $\omega_{\alpha} ( z ) $, $\omega_{\beta
} ( z ) $ and $ ( \omega_{\alpha}(z)+\omega_{\beta
}(z)-z ) ^{-2}$ in the rectangle $R_{\varepsilon}$, we have
%
\begin{equation}
\biggl\llvert \frac{m_{\mu_{\alpha}}^{\prime}(\omega_{\beta} (
z ) )+m_{\mu_{\beta}}^{\prime}(\omega_{\alpha} (
z ) )}{%
(\omega_{\alpha} ( z ) +\omega_{\beta} ( z )
-z)^{2}}%
-m_{\mu_{\alpha}}^{\prime}
\bigl(\omega_{\beta} ( z ) \bigr)m_{\mu
_{\beta}}^{\prime}\bigl(
\omega_{\alpha} ( z ) \bigr)\biggr\rrvert \geq c>0,
\end{equation}
everywhere in $R_{\varepsilon}$, provided that $\varepsilon$ is chosen
sufficiently small.

By using Lemma \ref{lemmaclosenessSttransforms}, we conclude that the
determinant
\[
\bigl\llvert \bigl[ m_{A}^{\prime} ( w_{2} )
+m_{B}^{\prime
} ( w_{1} ) \bigr] (
w_{1}+w_{2}-z ) ^{-2}-m_{A}^{\prime}
( w_{2} ) m_{B}^{\prime} ( w_{1} ) \bigr
\rrvert \geq c>0,
\]
where $w_{1}=\omega_{\alpha} ( z ) $, $w_{2}=\omega
_{\beta
} ( z ) $ and $z\in R_{\varepsilon}$.

It follows (with some additional help from Lemma \ref{lemmaclosenessSttransforms}), that the entries of the matrix $ [
F^{\prime} ] ^{-1}$ are bounded at $ ( \omega_{\alpha
},\omega
_{\beta} ) $ if $\operatorname{Im}z$ is sufficiently small [so that
$\operatorname{Im} ( w_{1}+w_{2}-z ) $ is bounded away from zero]. By compactness
of $R_{\varepsilon}$ and continuity of entries of the matrix $ [
F^{\prime} ] ^{-1}$, this shows that the operator norm of $ [
F^{\prime} ] ^{-1}$ is bounded at $ ( \omega_{\alpha
},\omega
_{\beta} ) $ uniformly for $z\in R_{\varepsilon}$.

By a similar argument, an application of Lemma \ref{lemmaclosenessSttransforms} shows that the operator norm of
$F^{\prime
\prime}$ is bounded for all $ ( w_{1},w_{2} ) $ in a fixed
neighborhood of $ ( \omega_{\alpha},\omega_{\beta} ) $,
and the
bound is uniform on $R_{\varepsilon}$. For example, we can compute
\[
\frac{\partial^{2}F_{1}}{ ( \partial w_{2} ) ^{2}}=2 ( w_{1}+w_{2}-z )
^{-3}+m_{A}^{\prime\prime} ( w_{2} ),
\]
and this is uniformly bounded in a certain neighborhood of
$w_{1}=\omega
_{\alpha} ( z ) $, $w_{2}=\omega_{\beta} ( z )
$ if $%
z\in R_{\varepsilon}$ and $\varepsilon$ is sufficiently small. The crucial
fact here is that the imaginary parts of $\omega_{\alpha} (
z ) $
and $\omega_{\beta} ( z ) $ are uniformly bounded away
from zero
for all $z\in R_{\varepsilon}$.

This shows that conditions (i) and (iii) of the Kantorovich theorem are
satisfied with some $C_{0}$ and $M_{0}$. Since $\llVert \Gamma
_{0}F ( x_{0} ) \rrVert \leq C_{0}\llVert  F (
x_{0} ) \rrVert $, we define $\delta_{0}:=C_{0}\llVert
F (
x_{0} ) \rrVert $ and note that $\delta_{0}=O (
r+s ) $.
By appropriate choice of $\overline{r}$ and $\overline{s}$, one can
make sure
that $h_{0}=C_{0}\delta_{0}M_{0}<1/2$ and, therefore, that conditions
(ii) and (iv) are satisfied. Moreover, one can make sure that $h_{0}$ is
arbitrarily small, and therefore that the neighborhood in the
conclusion of
the Kantorovich theorem has the form $\llVert  x-x_{0}\rrVert
\leq
\delta_{0}=O ( r+s ) $.

It follows by the Newton--Kantorovich theorem that there exists a
solution of
the equation $F ( w ) =0$ which satisfies the inequalities
\[
\bigl\llvert w_{1} ( z ) -\omega_{\alpha} ( z ) \bigr\rrvert
=O ( r+s )\quad\mbox{and}\quad\bigl\llvert w_{2} ( z ) -
\omega_{\beta} ( z ) \bigr\rrvert =O ( r+s ).
\]
The functions $\omega_{A} ( z ) $ and $\omega_{B} (
z ) $
defined by (\ref{omegadefinition}) satisfy equation $F ( w
) =0$,
and one can show that for every fixed $z$ they approach $\omega
_{\alpha
} ( z ) $ and $\omega_{\beta} ( z ) $ as
$N\rightarrow
\infty$. Hence, for sufficiently small $\overline{r}$ and $\overline{s}$
the solution of $F ( w ) =0$ found by the
Newton--Kantorovich method
coincide with the pair ($\omega_{A} ( z ),\omega_{B} (
z ) $) and we can conclude that
\[
\bigl\llvert \omega_{A} ( z ) -\omega_{\alpha} ( z ) \bigr
\rrvert =O ( r+s )\quad\mbox{and}\quad\bigl\llvert \omega_{B} ( z ) -
\omega_{\beta} ( z ) \bigr\rrvert =O ( r+s ).
\]

This completes the proof of Proposition \ref{propdiffsubordination} and
Theorem \ref{Thmdelocalization}.
\end{pf*}

\section{Local law for eigenvalues}\label{sec4}
\label{sectionlocallaw}

Let $\eta^{\ast}=M\eta$ and $I_{\eta^{\ast}}= [ x-\eta^{\ast
}+i\eta,x+\eta^{\ast}+i\eta ] $. Recall that
\[
s ( A,B ):=\max\bigl\{d_{L}(\mu_{A},
\mu_{\alpha}),d_{L}(\mu _{B},\mu_{\beta})
\bigr\}.
\]

\begin{proposition}
\label{propdiffStieltjesstoch}Assume that a pair of probability
measures $%
 ( \mu_{\alpha},\mu_{\beta} ) $ is smooth in a closed
interval $%
I$. Assume that $s ( A,B ) \leq\overline{s}$ where
$\overline{s}$
is a positive constant. Let $\eta=N^{-1/7}\log N$. Then for some
positive $c$
and $c_{1}$, and for every $\varepsilon>0$,
%
\begin{equation}
\quad\mathbb{P} \Bigl\{ \sup_{z\in I_{\eta^{\ast}}}\bigl\llvert m_{H} (
z ) -m_{\mu_{\alpha}\boxplus\mu_{\beta}} ( z ) \bigr\rrvert >\varepsilon+cs ( A,B ) \Bigr\} \leq
\exp \bigl( -c_{1} ( \log N ) ^{2} \bigr)
\end{equation}
for all sufficiently large $N$.
\end{proposition}

\begin{pf}
Proposition \ref{propdiffStieltjesstoch} is proved by combining Lemmas \ref{propdiffStieltjes} and \ref{lemmasupmlargedeviations} below.

\begin{lemma}
\label{propdiffStieltjes}Assume that a pair of probability measures $
 ( \mu_{\alpha},\mu_{\beta} ) $ is smooth in a closed
interval $%
I$. Let $r$ and $s$ be as defined in (\ref{definitionr}) and
(\ref{definitions}), respectively. Then for all sufficiently small $r$, $s$
and $\eta$, %
\[
\bigl\llvert \mathbb{E}m_{H} ( x+i\eta ) -m_{\mu_{\alpha
}\boxplus
\mu_{\beta}} ( x+i\eta
) \bigr\rrvert <O ( r+s ).
\]
\end{lemma}

Indeed, since $\mathbb{E}m_{H}= ( \omega_{A}+\omega_{B}-z
) ^{-1}$
and $m_{\mu_{\alpha}\boxplus\mu_{\beta}}= ( \omega_{\alpha
}+\omega_{\beta}-z ) ^{-1}$, therefore,
\[
\mathbb{E}m_{H}-m_{\mu_{\alpha}\boxplus\mu_{\beta}}=\frac{\omega
_{\alpha}+\omega_{\beta}-\omega_{A}-\omega_{B}}{ ( \omega
_{A}+\omega_{B}-z )  ( \omega_{\alpha}+\omega_{\beta
}-z )
}.
\]
The denominator is bounded away from zero for small $\eta$ by
Proposition %
\ref{propdiffsubordination} [$\operatorname{Im}\omega_{\alpha} (
x ) $
and $\operatorname{Im}\omega_{\beta} ( x ) $ are bounded away
from $0$
by the assumption of Proposition \ref{propdiffsubordination}, and
$\omega
_{A}$ and $\omega_{B}$ are close to $\omega_{\alpha}$ and $\omega
_{\beta
}$, respectively, by its conclusion]. The numerator can be estimated by
Proposition \ref{propdiffsubordination} as $O ( r+s ) $.

In \cite{kargin11b}, the following result was proved (as Corollary 6).

\begin{lemma}
\label{lemmasupmlargedeviations} For some positive $c$ and $c_{1}$ which
can depend on $M$, and for all $\delta>0$,
\[
\mathbb{P} \Bigl\{ \sup_{z\in I_{\eta^{\ast}}}\bigl\llvert m_{H} (
z ) -\mathbb{E}m_{H} ( z ) \bigr\rrvert >\delta \Bigr\} \leq \exp
\biggl( -\frac{c\delta^{2}\eta^{4}}{\llVert  B\rrVert
^{2}}%
N^{2} \biggr),
\]
provided that $N\geq c_{1} ( \sqrt{-\log ( \eta\delta
 ) }%
 ) / ( \eta^{2}\delta ) $.
\end{lemma}

Let us take $\delta=c\log N/ ( N\eta^{2} ) $ in Lemma \ref{lemmasupmlargedeviations}. Then $N\geq c_{1} ( \sqrt{-\log
 (
\eta\delta ) } ) /\break ( \eta^{2}\delta ) $
provided that
$\eta\geq N^{-1}\log N$. In particular, if $\eta=N^{-1/7}\log N$, then
Lemma \ref{lemmasupmlargedeviations} implies that
\[
\mathbb{P} \biggl\{ \sup_{z\in I_{\eta^{\ast}}}\bigl|m_{H}(z)-\mathbb
{E}m_{H}(z)\bigr|>%
\frac{1}{N^{5/7}\log N} \biggr\} \leq\exp \bigl( -c
( \log N ) ^{2} \bigr).
\]

In addition, the second part of Theorem \ref{PropRAestimate} and the
definition of $r$ imply that if $\eta\gg N^{-1/7}$, then for every $%
\varepsilon>0$ and all sufficiently large $N$, we have $\llvert
r (
z ) \rrvert <\varepsilon$. Hence, Lemma \ref{propdiffStieltjes}
implies that if $\eta\gg N^{-1/7}$, then
\[
\bigl\llvert \mathbb{E}m_{H} ( z ) -m_{\mu_{\alpha
}\boxplus\mu
_{\beta}} ( z ) \bigr
\rrvert <\varepsilon+cs ( A,B )
\]
for all sufficiently large $N$. Together, these statements imply the claim
of Proposition \ref{propdiffStieltjesstoch}.
\end{pf}

\begin{pf*}{Proof of Theorem \ref{Thmlocallaw}}
The proof is similar
to the
proof of Corollary 4.2 in \cite{erdosschleinyau09}. Let $\eta=cN^{-1/7}$,
and $c$ is sufficiently large, and let $\eta^{\ast}=M\eta$. Let
\begin{eqnarray*}
R ( \lambda ) &:=&\frac{1}{\pi}\int_{x-\eta^{\ast
}}^{x+\eta
^{\ast}}
\frac{\eta}{ ( x-\lambda ) ^{2}+\eta^{2}}\,dx
\\
&=&\frac{1}{\pi} \biggl( \arctan \biggl( \frac{x-\lambda}{\eta
}+M \biggr) -
\arctan \biggl( \frac{x-\lambda}{\eta}-M \biggr) \biggr).
\end{eqnarray*}

Then $R=1_{I^{\ast}}+T_{1}+T_{2}+T_{3}$, where $1_{I^{\ast}}$ is the
indicator function of the interval $I^{\ast}= [ x-\eta^{\ast
},x+\eta
^{\ast} ] $ and functions $T_{1}$, $T_{2}$ and $T_{3}$ satisfy the
following properties:
\begin{eqnarray*}
\llvert T_{1}\rrvert &\leq& c/\sqrt{M},\qquad\operatorname{supp}%
(T_{1})
\subset I_{1}= \bigl[ x-2\eta^{\ast},x+2\eta^{\ast}
\bigr],
\\
\llvert T_{2}\rrvert &\leq&1,\qquad\operatorname{supp}(T_{2})
\subset J_{1}\cup J_{2},
\end{eqnarray*}
where $J_{1}$ and $J_{2}$ are intervals of length $\sqrt{M}\eta$ with
midpoints at $x-\eta^{\ast}$ and $x+\eta^{\ast}$, respectively, and
\[
|T_{3}|\leq\frac{C\eta\eta^{\ast}}{ ( \lambda-x )
^{2}+ (
\eta^{\ast} ) ^{2}},\qquad\operatorname{supp}(T_{3})
\in I_{1}^{c}.
\]

Note that
\begin{eqnarray*}
\frac{\mathcal{N}_{\eta^{\ast}} ( x ) }{2\eta^{\ast}N} &=&\frac{%
1}{2\eta^{\ast}}\int1_{I^{\ast}} ( \lambda ) \mu
_{H_{N}} ( d\lambda )
\\
&=&\frac{1}{2\eta^{\ast}}\int R ( \lambda ) \mu _{H_{N}} ( d\lambda ) -
\frac{1}{2\eta^{\ast}}\int ( T_{1}+T_{2}+T_{3} )
\mu_{H_{N}} ( d\lambda ).
\end{eqnarray*}
The last integral can be estimated as follows:
\[
\frac{1}{2\eta^{\ast}}\int\llvert T_{1}+T_{2}+T_{3}
\rrvert \mu _{H_{N}} ( d\lambda ) \leq\frac{c}{\sqrt{M}}
\frac
{\mathcal{N}%
_{I_{1}}}{2\eta^{\ast}N}+\frac{\mathcal{N}_{J_{1}}+\mathcal
{N}_{J_{2}}}{%
2\eta^{\ast}N}+\frac{C\eta}{\eta^{\ast}}\rho_{\eta_{\ast
}} ( x ),
\]
where $\mathcal{N}_{I}$ denote the number of eigenvalues of $H_{N}$ in
interval $I$, and
\[
\rho_{\eta^{\ast}} ( x ):=\frac{1}{\pi}\operatorname{Im}%
m_{H_{N}}
\bigl( x+i\eta^{\ast} \bigr) =\frac{1}{\pi}\int\frac
{\eta
^{\ast}}{ ( x-\lambda ) ^{2}+ ( \eta^{\ast} )
^{2}}\mu
_{H_{N}} ( d\lambda ).
\]
Hence, by using the inequality $\mathcal{N}_{\eta} ( x )
\leq
CN\eta\rho_{\eta} ( x ) $, one obtains
%
\begin{eqnarray}\label{remainderestimate}
&& \frac{1}{2\eta^{\ast}}\int\llvert T_{1}+T_{2}+T_{3}
\rrvert \mu _{H_{N}} ( d\lambda )
\nonumber\\[-8pt]\\[-8pt]
&&\qquad \leq\frac{c}{\sqrt{M}} \bigl( \rho
_{2\eta
^{\ast}} ( x ) +\rho_{\sqrt{M}\eta} \bigl( x-\eta^{\ast
} \bigr)
+\rho_{\sqrt{M}\eta} \bigl( x+\eta^{\ast} \bigr) +\rho_{\eta
_{\ast
}} (
x ) \bigr). \nonumber
\end{eqnarray}

By the second path of Theorem \ref{PropRAestimate}, $\mathbb{E}%
m_{H_{N}} ( x+i\eta^{\ast} ) -m_{A_{N}} ( \omega
_{B} (
x+i\eta^{\ast} )  ) =O ( \frac{1}{N^{2}\eta
^{6}} )
=o ( 1 ) $. In addition, the assumption of smoothness and
formula (%
\ref{omegaBandomegabeta}) imply that $m_{A_{N}} ( \omega
_{B} (
x+i\eta^{\ast} )  ) $ is bounded for every $x$ in the
interval $%
I $. Hence, the integral in (\ref{remainderestimate}) is bounded by
$O (
M^{-1/2} ) $.

The main term can be written as
\begin{eqnarray*}
&& \frac{1}{2\eta^{\ast}}\int_{I^{\ast}}\frac{1}{\pi}\operatorname{Im}m_{\mu
_{\alpha}\boxplus\mu_{\beta}} ( x+i\eta ) \,dx
\\
&&\qquad {} +\frac
{1}{2\eta
^{\ast}}\int_{I^{\ast}}
\frac{1}{\pi}\operatorname{Im} \bigl( m_{H_{N}} ( x+i\eta )
-m_{\mu_{\alpha}\boxplus\mu_{\beta}} ( x+i\eta ) \bigr) \,dx.
\end{eqnarray*}
The first part converges to $\rho_{\mu_{\alpha}\boxplus\mu_{\beta
}} ( x ) $ because the assumption that $ ( \mu_{\alpha
},\mu
_{\beta} ) $ is smooth at $x$ implies that $\mu_{\alpha
}\boxplus\mu
_{\beta}$ has an analytic density in a neighborhood of $x$. For the second
term, we can use the estimate in Proposition \ref{propdiffStieltjesstoch}.
The assumption that $s ( A_{N},B_{N} ) =\max \{
d_{L}(\mu
_{A_{N}},\mu_{\alpha}),d_{L}(\mu_{B_{N}},\mu_{\beta}) \}
\rightarrow0$ implies that this term converges to $0$ in probability
as $%
N\rightarrow\infty$.
\end{pf*}

\section{Subordination and spikes}\label{sec5}
\label{sectionspikes}

The proof of Theorem \ref{theoremspikes} is based on the following result.

\begin{proposition}
\label{propositionformulaspikes}Let $H=A+UBU^{\ast}$ where $A$ is a
rank-one $N$-by-$N$ Hermitian matrix with the nonzero eigenvalue
$\theta$,
and $B$ is an $N$-by-$N$ Hermitian matrix with the empirical eigenvalue
distribution $\mu_{B}$. Then the expected Stieltjes transform of $H$
satisfies the following equation for every $z\in\mathbb{C}^{+}$:
%
\begin{equation}
\mathbb{E}m_{H} ( z ) =m_{B} ( z ) +\frac
{1}{N}
\frac{%
m_{B}^{\prime} ( z ) }{m_{B} ( z ) } \biggl( \frac{1}{%
\theta m_{B} ( z ) +1}-1 \biggr) +O_{\eta} \biggl(
\frac
{1}{N^{2}}%
 \biggr). \label{formulaspikes}
\end{equation}
\end{proposition}

Here $O_{\eta} ( N^{-2} ) $ denotes a function $f (
z ) $
such that $N^{2}\llvert  f ( z ) \rrvert \leq C
( \operatorname{Im}z ) ^{-k}$ for some $k>0$ and $C>0$.

\begin{pf*}{Proof of Proposition \ref{propositionformulaspikes}}
Note that
\[
m_{A} ( z ) =-\frac{1}{z}+\frac{1}{N} \biggl(
\frac
{1}{\theta-z}+%
\frac{1}{z} \biggr).
\]
From Theorem \ref{PropRAestimate}, we know that the following system
holds for $\mathbb{E}m_{H} ( z ) $, $\omega_{A} (
z ) $, $%
\omega_{B} ( z ) $:
%
\begin{eqnarray}\label{systemspikes}
\mathbb{E}m_{H} ( z ) &=&-\frac{1}{\omega_{B} (
z ) }+%
\frac{1}{N} \biggl( \frac{1}{\theta-\omega_{B} ( z )
}+\frac{1}{%
\omega_{B} ( z ) } \biggr)
+O_{\eta} \bigl( N^{-2} \bigr), \nonumber
\\
\mathbb{E}m_{H} ( z ) &=&m_{B} \bigl(
\omega_{A} ( z ) \bigr) +O_{\eta} \bigl( N^{-2}
\bigr)\quad\mbox{and}
\\
z &=&\omega_{A} ( z ) +\omega_{B} ( z ) +
\frac
{1}{%
\mathbb{E}m_{H} ( z ) }.
\nonumber
\end{eqnarray}
This system can be considered as a perturbation of the system
%
\begin{eqnarray} \label{systemspikes0}
\overline{m}_{H}(z)+\frac{1}{\overline{\omega}_{B} ( z )
} &=&0,\nonumber
\\
\overline{m}_{H} ( z ) -m_{B} \bigl( \overline{\omega
}_{A} ( z ) \bigr) &=&0\quad\mbox{and}
\\
\overline{\omega}_{A} ( z ) +\overline{\omega}_{B} ( z ) +
\frac{1}{\overline{m}_{H} ( z ) } &=&z.
\nonumber
\end{eqnarray}
The solution of the unperturbed system is $\overline{m}_{H} (
z )
=m_{B} ( z ) $, $\overline{\omega}_{A} ( z )
=z$ and $%
\overline{\omega}_{B} ( z ) =-1/m_{B} ( z ) $.
We compute
the derivative of the unperturbed system (\ref{systemspikes0}) with
respect to $(\overline{m}_{H},\overline{\omega}_{A},\overline
{\omega}_{B})$
at the solution and find
%
\begin{eqnarray}\label{formulaforJ}
J&=&\pmatrix{
1 & 0 & -\dfrac{1}{\overline{\omega}_{B}^{2} ( z ) }
\vspace*{5pt}\cr
1 & -m_{B}^{\prime} \bigl( \overline{\omega}_{A} ( z ) \bigr) & 0
\vspace*{5pt}\cr
-\dfrac{1}{\overline{m}_{H}^{2} ( z ) } & 1 & 1%
}
\nonumber\\[-5pt]\\[-5pt]
&=&
\pmatrix{
1 & 0 & -m_{B}^{2}( z )
\vspace*{5pt}\cr
1 & -m_{B}^{\prime} ( z ) & 0
\vspace*{5pt}\cr
-\dfrac{1}{\overline{m}_{B}^{2} ( z ) } & 1 & 1}.\nonumber
\end{eqnarray}
From (\ref{systemspikes}), the perturbation of the system is
\[
\biggl( \dfrac{1}{N} \biggl( \dfrac{1}{\theta-\overline{\omega
}_{B} (
z ) }+\dfrac{1}{\overline{\omega}_{B} ( z ) } \biggr)
+O_{\eta
} \bigl( N^{-2} \bigr),O_{\eta} \bigl(
N^{-2} \bigr),0 \biggr).
\]
Note that
\[
\dfrac{1}{N} \biggl( \dfrac{1}{\theta-\overline{\omega}_{B} (
z ) }+%
\dfrac{1}{\overline{\omega}_{B} ( z ) }
\biggr) =\dfrac
{m_{B} (
z ) }{N} \biggl( \dfrac{1}{\theta m_{B} ( z )
+1}-1 \biggr).
\]
Hence, the linearized system is
\[
J\pmatrix{\Delta\overline{m}_{H}
\vspace*{5pt}\cr
\Delta\overline{\omega}_{A}
\vspace*{5pt}\cr
\Delta\overline{\omega}_{B}}
=
\pmatrix{
\dfrac{m_{B} ( z ) }{N} \biggl(
\dfrac{1}{\theta m_{B} (
z ) +1}-1 \biggr) +O_{\eta} \bigl( N^{-2} \bigr)
\vspace*{5pt}\cr
O_{\eta} \bigl( N^{-2} \bigr)
\vspace*{5pt}\cr
0},
\]
where $\Delta\overline{m}_{H}$, $\Delta\overline{\omega}_{A}$ and
$%
\Delta\overline{\omega}_{B}$ denote the first-order changes in the
solution caused by perturbation. By using this linearization and the formula
(\ref{formulaforJ}) for the derivative~$J$, we can easily compute the
linear approximation for the solution of the perturbed system. In
particular,
\[
\mathbb{E}m_{H} ( z ) =m_{B} ( z ) +\dfrac
{1}{N}
\dfrac{%
m_{B}^{\prime} ( z ) }{m_{B} ( z ) } \biggl( \dfrac{1}{%
\theta m_{B} ( z ) +1}-1 \biggr) +O_{\eta} \biggl(
\dfrac
{1}{N^{2}}%
 \biggr).
\]
The contribution of the higher order terms is $O_{\eta} (
N^{-2} )$.
\end{pf*}

\begin{pf*}{Proof of Theorem \ref{theoremspikes}}
Proof of this theorem is
similar to the proof of Theorem 2.1 in \cite{capitainedonati-martinferal09}
and for this reason we will be concise. Let us start with the case when
$%
\rho_{\mu} ( \theta_{0} ) >L$. Then the first step is to show
that for large $N$ there are no eigenvalues of $H_{N}$ in
$S_{\varepsilon
}:= ( L+\varepsilon,\rho_{\mu} ( \theta_{0} )
-\varepsilon
 ) \cup ( \rho_{\mu} ( \theta_{0} )
+\varepsilon,\infty ) $. In order to do this, we note that for all sufficiently
large $N$, the first correction term in formula~(\ref{formulaspikes}),
\[
\mathcal{L}_{N} ( z ):=\frac{m_{B_{N}}^{\prime
}}{m_{B_{N}}} \biggl(
\frac{1}{\theta m_{B_{N}}+1}-1 \biggr),
\]
is the Stieltjes transform of a distribution $\Lambda_{B_{N}}$ with a
compact support which must be outside of $S_{\varepsilon}$.
Verification of
this fact can be done as in the proof of Proposition~4.5 in \cite{capitainedonati-martinferal09}.

Next, one can use the Stieltjes inversion formula, which holds for
distributions by the results of Tillmann in \cite{tillmann53}.
Applying it
to formula (\ref{formulaspikes}), one finds that for every $\varphi
\in
C_{c}^{\infty} ( \mathbb{R} ) $,
\begin{eqnarray*}
&& \mathbb{E} \bigl[ N^{-1}\operatorname{Tr} \bigl( \varphi ( H_{N} )
\bigr)  \bigr]
\\
&&\qquad  =\int\varphi \,d\mu_{B_{N}}+\frac{1}{N}
\Lambda_{B_{N}} ( \varphi ) -\frac{1}{\pi}\lim_{\eta\rightarrow0^{+}}
\operatorname{Im}\int_{%
\mathbb{R}}\varphi ( x ) f ( x+i\eta ) \,dx,
\end{eqnarray*}
where $f ( x ) $ denotes the error term in (\ref
{formulaspikes}), $%
f ( x ) =O_{\eta} ( N^{-2} ) $. The last term is
$O (
N^{-2} ) $ (see Section~6 in \cite{haagerupthorbjornsen05} or
the Appendix in \cite{capitainedonati-martin07}) and, therefore, we find that
\[
\mathbb{E} \bigl[ N^{-1}\operatorname{Tr} \bigl( \varphi ( H_{N} )
\bigr) %
 \bigr] =\int\varphi \,d\mu_{B_{N}}+\frac{1}{N}
\Lambda_{B_{N}} ( \varphi ) +O \bigl( N^{-2} \bigr).
\]
In particular, if the support of $\varphi$ is in $S_{\varepsilon}$, then
the first and the second terms are zero and $\mathbb{E} [
N^{-1}\operatorname{Tr} ( \varphi ( H_{N} )  )  ] =O (
N^{-2} )
$. If in addition $\varphi$ is nonnegative, then by the Markov inequality
\[
\mathbb{P} \biggl[ N^{-1}\operatorname{Tr} \bigl( \varphi ( H_{N}
) \bigr) >\frac{1}{2N} \biggr] <\frac{\mathbb{E} [ N^{-1}\operatorname{Tr} ( \varphi
 ( H_{N} )  )  ] }{2N}=O \biggl(
\frac
{1}{N} \biggr).
\]
By using a sequence of functions $\varphi$ that approximate the indicator
function of $S_{\varepsilon}$, it follows that
\[
\mathbb{P} [\mbox{there is an eigenvalue of }H_{N}\mbox{ in
}%
S_{\varepsilon} ] <\frac{c}{N}.
\]

The next step is to show that for sufficiently large $N$, there is exactly
one eigenvalue to the right of $\rho_{\mu} ( \theta_{0} )
-\varepsilon$. This can be done similarly to the corresponding result
(Theorem 4.5) in \cite{capitainedonati-martinferal09}. Namely, note
that $%
\rho_{\mu} ( \theta ) $ is an increasing function for
$\theta>
\theta_{0}-\varepsilon$. [This follows from the fact that $m_{\mu
} (
x ) $ is a decreasing function for $x>L$.] Hence, we can find an
interval $ [ \alpha,\beta ] $ in $ ( \theta
_{0}-\varepsilon,\theta_{0} ) $ that will map to an interval $ [ a,b
] $ in $%
 ( L,\rho_{\mu} ( \theta_{0} )  ) $, with
$a:=\rho
_{\mu} ( \alpha ) $ and $b:=\rho_{\mu} ( \beta
 ) $.
The claim is that if $\lambda_{1} ( H_{N} ) $ and $\lambda
_{2} ( H_{N} ) $ are the largest and the second largest eigenvalues
of $H_{N}$, then
\[
\mathbb{P} \bigl[ \lambda_{2} ( H_{N} ) <a\mbox{ and }
\lambda _{1} ( H_{N} ) >b \bigr] \rightarrow1
\]
as $N\rightarrow\infty$.

From interlacing inequalities for matrices, we immediately obtain that  $
\lambda_{2} ( H_{N} ) <a$ for sufficiently large $N$. In
order to
prove that $\lambda_{1} ( H_{N} ) >b$, we consider matrix $%
cA_{N}+B_{N}$.  By using Weyl's inequalities and the uniform bound on
norms of $B_{N}$ we obtain that the largest eigenvalue $\lambda
_{1} (
cA_{N}+B_{N} ) \geq c\theta-\delta$ for some positive $\delta
$. On
the other hand $\rho_{\mu} ( c\beta ) \sim c\beta$ for
large $%
c$. We conclude that
\[
\lambda_{1} ( \overline{c}A_{N}+B_{N} ) >
\rho_{\mu} ( \overline{c}\beta )
\]
for a sufficiently large $\overline{c}$. In addition, Weyl's inequalities
imply that $\lambda_{1} ( c_{1}A_{N}+B_{N} ) -\lambda
_{1} (
c_{2}A_{N}+B_{N} ) \leq\llvert  c_{1}-c_{2}\rrvert
\theta_{0}$.
Hence, if $c$ changes slowly, then the first eigenvalue of $cA_{N}+B_{N}$
changes slowly. By what we proved above, there are no eigenvalues of $%
cA_{N}+B_{N}$ in the interval $ ( L+\varepsilon,\rho_{\mu} (
c\theta_{0} ) -\varepsilon ) $ with large probability.
Since $%
\rho_{\mu} ( c\theta_{0} ) $ is an increasing function of $c$
for $c\geq1$, hence the length of this interval is always $\geq\rho
_{\mu
} ( \theta_{0} ) -L-2\varepsilon>\varepsilon^{\prime
}>0$. By
changing $c$ along a finite sequence $\overline{c}=c_{1}>c_{2}>\cdots
>c_{l}=1$ with $\llvert  c_{i}-c_{i+1}\rrvert \leq\varepsilon
^{\prime}/\theta_{0}$, we can ensure that $\lambda_{1} (
c_{i}A_{N}+B_{N} ) >\rho_{\mu} ( c_{i}\beta ) $ for
all $i$
with large probability. Hence, as $N$ grows, the probability that
$\lambda
_{1} ( H_{N} ) \geq\rho_{\mu} ( \beta ) >\rho
_{\mu
} ( \theta ) -\varepsilon$ approaches $1$. Together with
the fact
that with high probability the interval $ ( L+\varepsilon,\rho
_{\mu
} ( \theta ) -\varepsilon ) \cup ( \rho_{\mu
} (
\theta ) +\varepsilon,\infty ) $ contains no
eigenvalues, this
implies that $\lambda_{1}$ converges in probability to $\rho_{\mu
} (
\theta ) $ as $N\rightarrow\infty$.

Next, consider the case when $\rho_{\mu} ( \theta_{0} )
\leq L$.
Then we conclude (by the argument at the start of the proof) that for every
fixed $\varepsilon>0$ there are no eigenvalues of $H_{N}$ in $%
S_{\varepsilon}:= ( L+\varepsilon,\infty ) $ with high
probability for large $N$. On the other hand, by Weyl's inequalities $%
\lambda_{1} ( H_{N} ) \geq\lambda_{1} ( B_{N}
) $.
Since $\lambda_{1} ( B_{N} ) \rightarrow L$ in probability, we
conclude that $\lambda_{1} ( H_{N} ) \rightarrow L$ in
probability.
\end{pf*}

\begin{appendix}\label{app}
\section{A derivation of formula (\texorpdfstring{\protect\ref{13154845497151451}}{13})}\label{secB}

Let $G ( z ) \equiv G_{H} ( z ) = (
A+B-z )
^{-1}$, where $B=U\widetilde{B}U^{\ast}$ and $U$ is a uniformly
distributed unitary matrix. Let $B_{t}=e^{iXt}Be^{-iXt}$ where $X$ is
Hermitian and let $G_{t}= ( A+B_{t}-z ) ^{-1}$. Then
$\mathbb{E}%
_{U} ( dG_{t}/dt ) =0$ for every Hermitian matrix $X$. Let
us for
clarity omit the subscript $U$ in the expectations below and treat $A$ as
fixed. It is easy to compute that $\partial G_{uv}/\partial
B_{xy}=-G_{ux}G_{yv}$ and that $dB_{t}/dt=i [ X,B ] $. By
using the
chain rule, we calculate $dG_{t}/dt$ and infer that
\[
\mathbb{E} \bigl( ( G_{H} ) _{ua} ( BG_{H} )
_{bv} \bigr) =%
\mathbb{E} \bigl( ( G_{H}B )
_{ua} ( G_{H} ) _{bv} \bigr).
\]
Setting $u=a$ and summing over all $a$ gives the identity
\[
\mathbb{E} ( m_{H}BG_{H} ) =\mathbb{E} (
f_{B}G_{H} ).
\]

It follows that
\begin{eqnarray*}
\mathbb{E} ( m_{H}G_{H} ) &=&\mathbb{E} (
m_{H}G_{A}-m_{H}G_{A}BG_{H}
)
\\
&=&\mathbb{E} ( m_{H}G_{A}-G_{A}f_{B}G_{H}
),
\end{eqnarray*}
where we used the identity $G_{H} ( z ) =G_{A} (
z )
-G_{A} ( z ) BG_{H} ( z ) $ in the first line.

This can be written in the following equivalent form:
\begin{eqnarray*}
\mathbb{E}m_{H}\mathbb{E}G_{H} &=&(\mathbb{E}m_{H})G_{A}-(
\mathbb{E} 
f_{B})G_{A}\mathbb{E}G_{H}
\\
&&{} -\mathbb{E}\bigl[(m_{H}-\mathbb{E}m_{H})G_{H}
\bigr]-G_{A}\mathbb {E}\bigl[(f_{B}-\mathbb{E}%
f_{B})G_{H}
\bigr]
\\
&=&(\mathbb{E}m_{H})G_{A}-(\mathbb{E}f_{B})G_{A}
\mathbb {E}G_{H}+\mathbb{E}%
\Delta_{A},
\end{eqnarray*}
where
\[
\Delta_{A}=- ( m_{H}-\mathbb{E}m_{H} )
G_{H}-G_{A} ( f_{B}-%
\mathbb{E}f_{B} ) G_{H}.
\]
This expression can be further rewritten (after we multiply it by
$A-z$ and
rearrange terms) as
\[
\mathbb{E}m_{H} \biggl( A- \biggl( z-\frac{\mathbb{E}f_{B}}{\mathbb
{E}m_{H}}%
 \biggr) \biggr) \mathbb{E}G_{H}=\mathbb{E}m_{H}+ ( A-z )
\mathbb{E}%
\Delta_{A}.
\]
Let $z^{\prime}:=z-\mathbb{E}f_{B}/\mathbb{E}m_{H}$. Then
\[
\mathbb{E}m_{H}\mathbb{E}G_{H}=G_{A} \bigl(
z^{\prime} \bigr) \mathbb{E}%
m_{H}+ ( A-z )
G_{A} \bigl( z^{\prime} \bigr) \mathbb {E}\Delta_{A}.
\]
Divide the resulting expression by $\mathbb{E}m_{H}$. Then we obtain
\begin{eqnarray*}
\mathbb{E}G_{H} ( z ) &=&G_{A} \bigl( z^{\prime}
\bigr) +\frac{1}{%
\mathbb{E}m_{H}} \bigl( ( A-z ) G_{A} \bigl( z^{\prime
}
\bigr) \mathbb{E}\Delta_{A} \bigr)
\\
&=&G_{A} \bigl( z^{\prime} \bigr) +R_{A}.
\end{eqnarray*}
%

\section{Some helpful lemmas about expected~resolvent}\label{secC}

The following result is from \cite{bbcf12}.

\begin{lemma}
\label{PropEGdiagonal} Suppose that $U$ is a uniformly distributed random
unitary matrix. Then $\mathbb{E} [  ( A+UBU^{\ast} )
^{-1}%
 ] $ belongs to the algebra generated by the matrix $A$. In particular,
if $A$ is diagonal, then $\mathbb{E} [  ( A+UBU^{\ast
} ) ^{-1}%
 ] $ is diagonal.
\end{lemma}

\begin{pf}
If $V$ is an arbitrary unitary matrix that commutes
with $A$,
then
\begin{eqnarray*}
V\mathbb{E} \bigl[ \bigl( A+UBU^{\ast} \bigr) ^{-1} \bigr]
V^{\ast
} &=&%
\mathbb{E} \bigl[ \bigl( VAV^{\ast}+VUB (
VU ) ^{\ast
} \bigr) ^{-1}%
 \bigr]
\\
&=&\mathbb{E} \bigl[ \bigl( A+UBU^{\ast} \bigr) ^{-1} \bigr].
\end{eqnarray*}
Hence, $\mathbb{E} [  ( A+UBU^{\ast} ) ^{-1} ]
$ commutes
with $V$. Since von Neumann algebras are generated by their unitaries, we
conclude that $\mathbb{E} [  ( A+UBU^{\ast} )
^{-1} ] $
belongs to the bicommutant of $A$. By the basic theorem about von Neumann
algebras, this bicommutant coincides with the algebra generated by $A$.
\end{pf}

Similarly, one can prove the following result.

\begin{lemma}
\label{lemmaexpectationandtrace} Suppose that $U$ is a uniformly
distributed random unitary matrix. Then
\[
\mathbb{E} \bigl[ UBU^{\ast} \bigr] = \biggl( \frac{1}{N}\operatorname{Tr} ( B ) \biggr) I_{N}.
\]
\end{lemma}

\begin{lemma}
\label{lemmaconvcombeigs}Let $A_{j}$, $j=1,\ldots,m$, be a family of
normal (finite-dimen\-sional) operators. Suppose that the eigenvalues of
all $%
A_{j}$ are contained in a closed disc $D\subset\mathbb{C}$, and let
$H=\sum
p_{j}A_{j}$ be a convex combination of $A_{j}$. Then all eigenvalues of $H$
are contained in $D$.
\end{lemma}

\begin{pf}
By subtracting a multiple of the identity operator
from all $%
A_{j}$, we can reduce the problem to the case when disc $D$ has its center
at $0$. Assume that this is indeed the case. Let $R$ be the radius of $D$.
Since the operators are normal, their norms are equal to the maximum of the
absolute values of eigenvalues. Hence, $\llVert  A_{j}\rrVert
\leq R$.
Hence, $\llVert  H\rrVert \leq\sum p_{j}\llVert
A_{j}\rrVert
\leq R$. It follows that all eigenvalues of $H$ have absolute value
$\leq R$.
\end{pf}

\section{Estimates of the resolvent entries, the Stieltjes transform and related quantities}\label{secD}

In this section, we assume that $G ( z ) =(A+UBU^{\ast}-z)^{-1}$,
where $A$ and $B$ are $N$-by-$N$ Hermitian matrices and $U$ is a random
Haar-distributed unitary matrix.

\begin{lemma}
\label{lemmadiffGiiEGii}Let $z=E+i\eta$ where $\eta>0$. Then, for a
numeric $c>0$ and every $\delta>0$:
\begin{longlist}[(ii)]
\item[(i)]
%
\begin{eqnarray}
\mathbb{P} \bigl\{ \bigl\llvert G_{ij}(z)-\mathbb{E}G_{ij} (
z ) \bigr\rrvert >\delta \bigr\} &\leq& \exp \biggl( -\frac{c\delta
^{2}\eta^{4}}{%
\llVert  B\rrVert ^{2}}N \biggr)
\quad\mbox{and}
\nonumber\\[-8pt]\\[-8pt]
\operatorname{Var} \bigl( G_{ij} ( z ) \bigr) &\leq&
\frac{\llVert  B\rrVert
^{2}}{c\eta
^{4}N};\nonumber
\end{eqnarray}

\item[(ii)] Let $h:=N^{-1}\operatorname{Tr} ( FG ) $, where $F$ does not depend
on $U$. Then
%
\begin{eqnarray}
\mathbb{P} \bigl\{ \bigl\llvert h ( z ) -\mathbb{E}h ( z ) \bigr\rrvert >\delta
\bigr\} &\leq& \exp \biggl( -\frac{c\delta
^{2}\eta^{4}}{%
\llVert  F\rrVert ^{2}\llVert  B\rrVert
^{2}}N^{2} \biggr)\quad
\mbox{and}
\nonumber\\[-8pt]\\[-8pt]
\operatorname{Var} \bigl( h ( z ) \bigr) &\leq&\frac{\llVert
F\rrVert ^{2}\llVert  B\rrVert ^{2}}{c\eta^{4}N^{2}}.\nonumber
\end{eqnarray}
\end{longlist}
\end{lemma}

\begin{rem*}
By applying the second part of the lemma to $h=I$, $A-z$
and $ ( A-z ) ^{-1}$, we can compute probabilities of deviations
and variances for $m ( z ):=N^{-1}\operatorname{Tr}G (
z ) $, $%
f_{B} ( z ):=N^{-1}\operatorname{Tr} ( BG ( z )
 )
=1-N^{-1}\operatorname{Tr} (  ( A-z ) G ( z )
 ) $ and
$h_{A} ( z ) =N^{-1}\operatorname{Tr} (  ( A-z )
^{-1}G ( z )  ) $, respectively. In particular,
%
\begin{eqnarray}
\mathbb{P} \bigl\{ \bigl\llvert m ( z ) -\mathbb{E}m ( z ) \bigr\rrvert >\delta
\bigr\} &\leq& \exp \biggl( -\frac{c\delta
^{2}\eta^{4}}{%
\llVert  B\rrVert ^{2}}N^{2} \biggr)\quad
\mbox{and}
\nonumber\\[-8pt]\\[-8pt]
\operatorname{Var} \bigl( m ( z ) \bigr) &\leq&\frac{\llVert  B\rrVert
^{2}}{c\eta
^{4}N^{2}};\nonumber
\\
\mathbb{P} \bigl\{ \bigl\llvert f_{B} ( z ) -\mathbb
{E}f_{B} ( z ) \bigr\rrvert >\delta \bigr\} &\leq& \exp \biggl[ -
\frac
{c\delta
^{2}\eta^{4}}{\llVert  A-z\rrVert ^{2}\llVert  B\rrVert ^{2}}%
N^{2} \biggr]
\quad\mbox{and}
\nonumber\\[-8pt]\\[-8pt]
\operatorname{Var} \bigl( f_{B} ( z ) \bigr) &\leq&\frac{\llVert  A-z\rrVert ^{2}\llVert  B\rrVert
^{2}}{c\eta
^{4}N^{2}};\nonumber
\\
\mathbb{P} \bigl\{ \bigl\llvert h_{A} ( z ) -\mathbb
{E}h_{A} ( z ) \bigr\rrvert >\delta \bigr\} &\leq&\exp \biggl[ -
\frac
{c\delta
^{2}\eta^{6}}{\llVert  B\rrVert ^{2}}N^{2} \biggr]
\quad\mbox {and}
\nonumber\\[-8pt]\\[-8pt]
\operatorname{Var} \bigl( h_{A} ( z ) \bigr) &\leq&
\frac{\llVert B\rrVert ^{2}}{c\eta^{6}N^{2}}.\nonumber
\end{eqnarray}
\end{rem*}

\begin{pf*}{Proof of Lemma \ref{lemmadiffGiiEGii}}
(i) In a small neighborhood of identity matrix, all unitary
matrices can be written as $U=e^{iX}$ where $X$ is Hermitian. Then $G_{H}$
can be thought of as a function of $X$ and we can compute its derivative as
follows [let $\widetilde{B}$ denote $UBU^{\ast}$, $B ( X
) =e^{iX}%
\widetilde{B}e^{-iX}$ and $G_{H} ( z,X ) = ( A+B (
X ) -z ) ^{-1}$]:
\begin{eqnarray*}
\bigl\llvert d_{X}G_{H} ( z,X ) \bigr\rrvert &=&\biggl
\llvert \sum_{x,y}%
\frac{\partial G_{H} ( z ) }{\partial\widetilde{B}_{xy}}%
\,d_{X}B_{xy}
( X ) \biggr\rrvert
\\
&=&\biggl\llvert \sum_{x,y}\frac{\partial G_{H} ( z )
}{\partial
\widetilde{B}_{xy}} [ X,
\widetilde{B} ] _{xy}\biggr\rrvert
\\
&=&\biggl\llvert \sum_{x,y} \biggl[
\frac{\partial G_{H} ( z
) }{%
\partial\widetilde{B}_{xy}},\widetilde{B} \biggr] X_{xy}\biggr\rrvert,
\end{eqnarray*}
where we used the fact that $ d_{X}B ( X ) \rrvert _{X=0}=
[ X,\widetilde{B} ] \equiv X\widetilde{B}-\widetilde{B}X$.

We compute
\[
\frac{\partial G_{ij}}{\partial\widetilde{B}_{xy}}=-G_{ix}G_{yj}.
\]
Therefore,
\[
\biggl\llVert \frac{\partial G_{ij}}{\partial\widetilde{B}_{xy}}\biggr\rrVert _{2}=\sqrt{\sum
_{x,y}\llvert G_{xi}\rrvert ^{2}
\llvert G_{yj}\rrvert ^{2}}=\sqrt{\llVert Ge_{i}
\rrVert ^{2}\llVert Ge_{j}\rrVert ^{2}}\leq\llVert
G\rrVert ^{2}\leq\frac
{1}{\eta
^{2}},
\]
where $\llVert  M \rrVert _{2}:=\operatorname{Tr} ( M^{\ast
}M ) $
is the Frobenius norm of matrix $M$.

If $\llVert  X\rrVert _{2}=1$, then it follows that
\begin{eqnarray*}
\bigl\llvert d_{X}G_{ij} ( z,X ) \bigr\rrvert &\leq&\biggl
\llVert \biggl[ \frac{\partial G_{ij} ( z ) }{\partial\widetilde{B}_{xy}}, 
\widetilde{B} \biggr] \biggr
\rrVert _{2}
\\
&\leq&2\biggl\llVert \frac{\partial G_{ij} ( z ) }{\partial
\widetilde{B}_{xy}}\biggr\rrVert _{2}\llVert B
\rrVert
\\
&\leq&\frac{2\llVert  B\rrVert }{\eta^{2}}.
\end{eqnarray*}
In the second line, we used the fact that $\llVert  AB\rrVert
_{2}\leq
\llVert  A\rrVert _{2}\llVert  B\rrVert $. (See
Exercise 20 on
page~313 in Section~5.6 of \cite{hornjohnson85}.)

Next, we note that the Ricci's curvature of $SU ( N ) $ is $(N/2)I$
with respect to the metric induced by $\llVert \cdot\rrVert _{2}$
norm on $X$. By Gromov's theorem, if $g\dvtx  ( \mathcal{SU} (
N ),\llVert  ds\rrVert _{2} ) \rightarrow\mathbb{R}$ is an
$L$-Lipschitz function and if $\mathbb{E}g=0$, then $P \{ \llvert
g\rrvert >\delta \} \leq\exp ( -cN\delta
^{2}/L^{2} ) $
for every $\delta>0$ and some numeric $c>0$. For details of the argument,
the reader can consult Section~4.4.2 in \cite{andersonguionnetzeitouni10},
especially Theorem 4.4.7. We apply this theorem to a complex-valued function
but the proof is the same except for some minor changes.

For variance, we note that for every positive random variable $X$, it is
true that $\mathbb{E}X=\int_{0}^{\infty} ( 1-\mathcal
{F}_{X} (
t )  ) \,dt$, where $\mathcal{F}_{X} ( t ) $ is cumulative
distribution function of $X$. We can apply this to the random variables
$%
 ( \operatorname{Im} ( G_{ij}-\mathbb{E}G_{ij} )  )
^{2}$ and $%
 ( \operatorname{Re} ( G_{ij}-\mathbb{E}G_{ij} )  )
^{2}$, and
find that the expectation of both expression is smaller than $\frac{%
\llVert  B\rrVert ^{2}}{c\eta^{4}N}$. Hence,
\[
\operatorname{Var} \bigl( G_{ij} ( z ) \bigr) \equiv\mathbb{E} 
 \bigl( ( G_{ij}-\mathbb{E}G_{ij} ) ( \overline
{G_{ij}-\mathbb{%
E}G_{ij}} ) \bigr) \leq
\frac{\llVert  B\rrVert
^{2}}{c\eta^{4}N%
}
\]
with a possibly different constant.

(ii) The\vspace*{2pt} proof is similar and boils down to showing that if $h (
z,X ):=N^{-1}\operatorname{Tr} ( FG_{H} ( z,X )
 ) $ and
if $\llVert  X\rrVert _{2}=1$, then
\begin{eqnarray*}
\bigl\llvert d_{X}h_{A} ( z,X ) \bigr\rrvert &=&\biggl
\llvert \frac{1}{N}%
\sum_{x,y} \bigl(
[ GFG,\widetilde{B} ] \bigr) _{yx}X_{xy}\biggr\rrvert
\\
&\leq&\frac{2}{\sqrt{N}}\llVert GFG\rrVert \llVert B\rrVert
\\
&\leq&\frac{2\llVert  F\rrVert \llVert  B\rrVert
}{\eta^{2}%
\sqrt{N}}.
\end{eqnarray*}\upqed
\end{pf*}

\begin{lemma}
\label{lemmaimmzlowerbound} Assume that $\max \{ \llVert
A\rrVert,\llVert  B\rrVert  \} \leq K$, $\operatorname{Im}z=\eta
>0$ and $\llvert  z\rrvert \leq R$. We have $ ( \mathbb
{E}%
m_{H} ( z )  ) ^{-1}\leq c^{\prime}/\eta$, where
$c^{\prime
} $ depends only on $K$ and $R$.
\end{lemma}

\begin{pf}
We have
\[
\operatorname{Im}\mathbb{E} \biggl[ \frac{1}{N}\operatorname{Tr}G_{H} ( x+i
\eta ) \biggr] =\mathbb{E} \biggl[ \frac{\eta}{N}\operatorname{Tr} \bigl[ \bigl( (
H-xI_{N} ) ^{2}+\eta^{2}I_{N} \bigr)
^{-1} \bigr] \biggr].
\]
Since all eigenvalues of $ ( H-xI_{N} ) ^{2}+\eta^{2}I_{N}$
are $%
\leq (  ( K+R ) ^{2}+R^{2} ) $, hence all
eigenvalues of $%
 (  ( H-xI_{N} ) ^{2}+\eta^{2}I_{N} ) ^{-1}$ are
${\geq} (  ( K+R ) ^{2}+R^{2} ) ^{-1}$ and, therefore,
\[
\mathbb{E} \biggl[ \frac{\eta}{N}\operatorname{Tr} \bigl[ ( H-xI_{N} )
^{2}+\eta^{2}I_{N} \bigr] ^{-1} \biggr]
\geq c\eta,
\]
which implies the claim of the lemma.
\end{pf}
\end{appendix}




\printaddresses

\begin{thebibliography}{44}

\bibitem{andersonguionnetzeitouni10}
\begin{bbook}[mr]
\bauthor{\bsnm{Anderson},~\bfnm{Greg~W.}\binits{G.~W.}},
\bauthor{\bsnm{Guionnet},~\bfnm{Alice}\binits{A.}} \AND
\bauthor{\bsnm{Zeitouni},~\bfnm{Ofer}\binits{O.}}
(\byear{2010}).
\btitle{An Introduction to Random Matrices}.
\bseries{Cambridge Studies in Advanced Mathematics}
\bvolume{118}.
\bpublisher{Cambridge Univ. Press},
\blocation{Cambridge}.
\bid{mr={2760897}}
\bptnote{check year}%
\end{bbook}
\bptok{imsref}%
\endbibitem

\bibitem{baiyao07}
\begin{barticle}[mr]
\bauthor{\bsnm{Bai},~\bfnm{Zhidong}\binits{Z.}} \AND
\bauthor{\bsnm{Yao},~\bfnm{Jian-feng}\binits{J.-f.}}
(\byear{2008}).
\btitle{Central limit theorems for eigenvalues in a spiked population model}.
\bjournal{Ann. Inst. Henri Poincar\'e Probab. Stat.}
\bvolume{44}
\bpages{447--474}.
\bid{doi={10.1214/07-AIHP118}, issn={0246-0203}, mr={2451053}}
\bptnote{check year}%
\end{barticle}
\bptok{imsref}%
\endbibitem

\bibitem{baikbenarouspeche05}
\begin{barticle}[mr]
\bauthor{\bsnm{Baik},~\bfnm{Jinho}\binits{J.}},
\bauthor{\bsnm{Ben Arous},~\bfnm{G{\'e}rard}\binits{G.}} \AND
\bauthor{\bsnm{P{\'e}ch{\'e}},~\bfnm{Sandrine}\binits{S.}}
(\byear{2005}).
\btitle{Phase transition of the largest eigenvalue for nonnull complex sample covariance matrices}.
\bjournal{Ann. Probab.}
\bvolume{33}
\bpages{1643--1697}.
\bid{doi={10.1214/009117905000000233}, issn={0091-1798}, mr={2165575}}
\end{barticle}
\bptok{imsref}%
\endbibitem

\bibitem{baiksilverstein06}
\begin{barticle}[mr]
\bauthor{\bsnm{Baik},~\bfnm{Jinho}\binits{J.}} \AND
\bauthor{\bsnm{Silverstein},~\bfnm{Jack~W.}\binits{J.~W.}}
(\byear{2006}).
\btitle{Eigenvalues of large sample covariance matrices of spiked population models}.
\bjournal{J. Multivariate Anal.}
\bvolume{97}
\bpages{1382--1408}.
\bid{doi={10.1016/j.jmva.2005.08.003}, issn={0047-259X}, mr={2279680}}
\end{barticle}
\bptok{imsref}%
\endbibitem

\bibitem{belinschi08}
\begin{barticle}[mr]
\bauthor{\bsnm{Belinschi},~\bfnm{Serban~Teodor}\binits{S.~T.}}
(\byear{2008}).
\btitle{The {L}ebesgue decomposition of the free additive convolution of two probability distributions}.
\bjournal{Probab. Theory Related Fields}
\bvolume{142}
\bpages{125--150}.
\bid{doi={10.1007/s00440-007-0100-3}, issn={0178-8051}, mr={2413268}}
\end{barticle}
\bptok{imsref}%
\endbibitem

\bibitem{belinschibercovici07}
\begin{barticle}[mr]
\bauthor{\bsnm{Belinschi},~\bfnm{S.~T.}\binits{S.~T.}} \AND
\bauthor{\bsnm{Bercovici},~\bfnm{H.}\binits{H.}}
(\byear{2007}).
\btitle{A new approach to subordination results in free probability}.
\bjournal{J. Anal. Math.}
\bvolume{101}
\bpages{357--365}.
\bid{doi={10.1007/s11854-007-0013-1}, issn={0021-7670}, mr={2346550}}
\end{barticle}
\bptok{imsref}%
\endbibitem

\bibitem{bbcf12}
\begin{bmisc}[auto:STB|2014/02/12|14:17:21]
\bauthor{\bsnm{Belinschi},~\bfnm{S.~T.}\binits{S.~T.}},
\bauthor{\bsnm{Bercovici},~\bfnm{H.}\binits{H.}},
\bauthor{\bsnm{Capitaine},~\bfnm{M.}\binits{M.}} \AND
\bauthor{\bsnm{F{\'e}vrier},~\bfnm{M.}\binits{M.}}
(\byear{2012}).
\bhowpublished{Outliers in the spectrum of large deformed unitarily invariant models.
Available at \arxivurl{arXiv:1207.5443}.}
\end{bmisc}
\bptok{imsref}%
\endbibitem

\bibitem{BenaychGeorgesGuionnetMaida11}
\begin{barticle}[mr]
\bauthor{\bsnm{Benaych-Georges},~\bfnm{F.}\binits{F.}},
\bauthor{\bsnm{Guionnet},~\bfnm{A.}\binits{A.}} \AND
\bauthor{\bsnm{Maida},~\bfnm{M.}\binits{M.}}
(\byear{2011}).
\btitle{Fluctuations of the extreme eigenvalues of finite rank deformations of random matrices}.
\bjournal{Electron. J. Probab.}
\bvolume{16}
\bpages{1621--1662}.
\bid{doi={10.1214/EJP.v16-929}, issn={1083-6489}, mr={2835249}}
\end{barticle}
\bptok{imsref}%
\endbibitem

\bibitem{BenaychGeorgesGuionnetMaida11a}
\begin{barticle}[mr]
\bauthor{\bsnm{Benaych-Georges},~\bfnm{F.}\binits{F.}},
\bauthor{\bsnm{Guionnet},~\bfnm{A.}\binits{A.}} \AND
\bauthor{\bsnm{Maida},~\bfnm{M.}\binits{M.}}
(\byear{2012}).
\btitle{Large deviations of the extreme eigenvalues of random deformations of matrices}.
\bjournal{Probab. Theory Related Fields}
\bvolume{154}
\bpages{703--751}.
\bid{doi={10.1007/s00440-011-0382-3}, issn={0178-8051}, mr={3000560}}
\end{barticle}
\bptok{imsref}%
\endbibitem

\bibitem{BenaychGeorgesNadakuditi11}
\begin{barticle}[mr]
\bauthor{\bsnm{Benaych-Georges},~\bfnm{Florent}\binits{F.}} \AND
\bauthor{\bsnm{Nadakuditi},~\bfnm{Raj~Rao}\binits{R.~R.}}
(\byear{2011}).
\btitle{The eigenvalues and eigenvectors of finite, low rank perturbations of large random matrices}.
\bjournal{Adv. Math.}
\bvolume{227}
\bpages{494--521}.
\bid{doi={10.1016/j.aim.2011.02.007}, issn={0001-8708}, mr={2782201}}
\end{barticle}
\bptok{imsref}%
\endbibitem

\bibitem{BenaychGeorgesNadakuditi12}
\begin{barticle}[mr]
\bauthor{\bsnm{Benaych-Georges},~\bfnm{Florent}\binits{F.}} \AND
\bauthor{\bsnm{Nadakuditi},~\bfnm{Raj~Rao}\binits{R.~R.}}
(\byear{2012}).
\btitle{The singular values and vectors of low rank perturbations of large rectangular random matrices}.
\bjournal{J. Multivariate Anal.}
\bvolume{111}
\bpages{120--135}.
\bid{doi={10.1016/j.jmva.2012.04.019}, issn={0047-259X}, mr={2944410}}
\end{barticle}
\bptok{imsref}%
\endbibitem

\bibitem{biane98b}
\begin{barticle}[mr]
\bauthor{\bsnm{Biane},~\bfnm{Philippe}\binits{P.}}
(\byear{1998}).
\btitle{Processes with free increments}.
\bjournal{Math. Z.}
\bvolume{227}
\bpages{143--174}.
\bid{doi={10.1007/PL00004363}, issn={0025-5874}, mr={1605393}}
\end{barticle}
\bptok{imsref}%
\endbibitem

\bibitem{bloemendalvirag11a}
\begin{bmisc}[mr]
\bauthor{\bsnm{Bloemendal},~\bfnm{Alex}\binits{A.}} \AND
\bauthor{\bsnm{Vir{\'a}g},~\bfnm{B{\'a}lint}\binits{B.}}
(\byear{2011}).
\bhowpublished{Limits of spiked random matrices {II}.
Available at \arxivurl{arXiv:1109.3704}.}
\end{bmisc}
\bptok{imsref}%
\endbibitem

\bibitem{bloemendalvirag11}
\begin{barticle}[mr]
\bauthor{\bsnm{Bloemendal},~\bfnm{Alex}\binits{A.}} \AND
\bauthor{\bsnm{Vir{\'a}g},~\bfnm{B{\'a}lint}\binits{B.}}
(\byear{2013}).
\btitle{Limits of spiked random matrices {I}}.
\bjournal{Probab. Theory Related Fields}
\bvolume{156}
\bpages{795--825}.
\bid{doi={10.1007/s00440-012-0443-2}, issn={0178-8051}, mr={3078286}}
\end{barticle}
\bptok{imsref}%
\endbibitem

\bibitem{capitaine12}
\begin{barticle}[mr]
\bauthor{\bsnm{Capitaine},~\bfnm{M.}\binits{M.}}
(\byear{2013}).
\btitle{Additive/multiplicative free subordination property and limiting eigenvectors of spiked additive deformations of {W}igner matrices and spiked sample covariance matrices}.
\bjournal{J. Theoret. Probab.}
\bvolume{26}
\bpages{595--648}.
\bid{doi={10.1007/s10959-012-0416-5}, issn={0894-9840}, mr={3090543}}
\end{barticle}
\bptok{imsref}%
\endbibitem

\bibitem{capitainedonati-martin07}
\begin{barticle}[mr]
\bauthor{\bsnm{Capitaine},~\bfnm{M.}\binits{M.}} \AND
\bauthor{\bsnm{Donati-Martin},~\bfnm{C.}\binits{C.}}
(\byear{2007}).
\btitle{Strong asymptotic freeness for {W}igner and {W}ishart matrices}.
\bjournal{Indiana Univ. Math. J.}
\bvolume{56}
\bpages{767--803}.
\bid{doi={10.1512/iumj.2007.56.2886}, issn={0022-2518}, mr={2317545}}
\end{barticle}
\bptok{imsref}%
\endbibitem

\bibitem{capitainedonati-martinferal09}
\begin{barticle}[mr]
\bauthor{\bsnm{Capitaine},~\bfnm{Mireille}\binits{M.}},
\bauthor{\bsnm{Donati-Martin},~\bfnm{Catherine}\binits{C.}} \AND
\bauthor{\bsnm{F{\'e}ral},~\bfnm{Delphine}\binits{D.}}
(\byear{2009}).
\btitle{The largest eigenvalues of finite rank deformation of large {W}igner matrices: Convergence and nonuniversality of the fluctuations}.
\bjournal{Ann. Probab.}
\bvolume{37}
\bpages{1--47}.
\bid{doi={10.1214/08-AOP394}, issn={0091-1798}, mr={2489158}}
\end{barticle}
\bptok{imsref}%
\endbibitem

\bibitem{capitainedonati-martinferalfevrier11}
\begin{barticle}[mr]
\bauthor{\bsnm{Capitaine},~\bfnm{M.}\binits{M.}},
\bauthor{\bsnm{Donati-Martin},~\bfnm{C.}\binits{C.}},
\bauthor{\bsnm{F{\'e}ral},~\bfnm{D.}\binits{D.}} \AND
\bauthor{\bsnm{F{\'e}vrier},~\bfnm{M.}\binits{M.}}
(\byear{2011}).
\btitle{Free convolution with a semicircular distribution and eigenvalues of spiked deformations of {W}igner matrices}.
\bjournal{Electron. J. Probab.}
\bvolume{16}
\bpages{1750--1792}.
\bid{doi={10.1214/EJP.v16-934}, issn={1083-6489}, mr={2835253}}
\end{barticle}
\bptok{imsref}%
\endbibitem

\bibitem{chistyakovgotze11}
\begin{barticle}[mr]
\bauthor{\bsnm{Chistyakov},~\bfnm{Gennadii~P.}\binits{G.~P.}} \AND
\bauthor{\bsnm{G{\"o}tze},~\bfnm{Friedrich}\binits{F.}}
(\byear{2011}).
\btitle{The arithmetic of distributions in free probability theory}.
\bjournal{Cent. Eur. J. Math.}
\bvolume{9}
\bpages{997--1050}.
\bid{doi={10.2478/s11533-011-0049-4}, issn={1895-1074}, mr={2824443}}
\end{barticle}
\bptok{imsref}%
\endbibitem

\bibitem{erdos10}
\begin{barticle}[mr]
\bauthor{\bsnm{Erd{\H{o}}s},~\bfnm{L.}\binits{L.}}
(\byear{2011}).
\btitle{Universality of {W}igner random matrices: A survey of recent results}.
\bjournal{Uspekhi Mat. Nauk}
\bvolume{66}
\bpages{67--198}.
\bid{doi={10.1070/RM2011v066n03ABEH004749}, issn={0042-1316}, mr={2859190}}
\end{barticle}
\bptok{imsref}%
\endbibitem

\bibitem{erdosknowles11}
\begin{barticle}[mr]
\bauthor{\bsnm{Erd{\H{o}}s},~\bfnm{L{\'a}szl{\'o}}\binits{L.}} \AND
\bauthor{\bsnm{Knowles},~\bfnm{Antti}\binits{A.}}
(\byear{2011}).
\btitle{Quantum diffusion and eigenfunction delocalization in a random band matrix model}.
\bjournal{Comm. Math. Phys.}
\bvolume{303}
\bpages{509--554}.
\bid{doi={10.1007/s00220-011-1204-2}, issn={0010-3616}, mr={2782623}}
\end{barticle}
\bptok{imsref}%
\endbibitem

\bibitem{ekyy12}
\begin{barticle}[mr]
\bauthor{\bsnm{Erd{\H{o}}s},~\bfnm{L{\'a}szl{\'o}}\binits{L.}},
\bauthor{\bsnm{Knowles},~\bfnm{Antti}\binits{A.}},
\bauthor{\bsnm{Yau},~\bfnm{Horng-Tzer}\binits{H.-T.}} \AND
\bauthor{\bsnm{Yin},~\bfnm{Jun}\binits{J.}}
(\byear{2013}).
\btitle{The local semicircle law for a general class of random matrices}.
\bjournal{Electron. J. Probab.}
\bvolume{18}
\bpages{no. 59, 58}.
\bid{doi={10.1214/EJP.v18-2473}, issn={1083-6489}, mr={3068390}}
\end{barticle}
\bptok{imsref}%
\endbibitem

\bibitem{erdosschleinyau09}
\begin{barticle}[mr]
\bauthor{\bsnm{Erd{\H{o}}s},~\bfnm{L{\'a}szl{\'o}}\binits{L.}},
\bauthor{\bsnm{Schlein},~\bfnm{Benjamin}\binits{B.}} \AND
\bauthor{\bsnm{Yau},~\bfnm{Horng-Tzer}\binits{H.-T.}}
(\byear{2009}).
\btitle{Semicircle law on short scales and delocalization of eigenvectors for {W}igner random matrices}.
\bjournal{Ann. Probab.}
\bvolume{37}
\bpages{815--852}.
\bid{doi={10.1214/08-AOP421}, issn={0091-1798}, mr={2537522}}
\end{barticle}
\bptok{imsref}%
\endbibitem

\bibitem{feralpeche07}
\begin{barticle}[mr]
\bauthor{\bsnm{F{\'e}ral},~\bfnm{Delphine}\binits{D.}} \AND
\bauthor{\bsnm{P{\'e}ch{\'e}},~\bfnm{Sandrine}\binits{S.}}
(\byear{2007}).
\btitle{The largest eigenvalue of rank one deformation of large {W}igner matrices}.
\bjournal{Comm. Math. Phys.}
\bvolume{272}
\bpages{185--228}.
\bid{doi={10.1007/s00220-007-0209-3}, issn={0010-3616}, mr={2291807}}
\end{barticle}
\bptok{imsref}%
\endbibitem

\bibitem{haagerupthorbjornsen05}
\begin{barticle}[mr]
\bauthor{\bsnm{Haagerup},~\bfnm{Uffe}\binits{U.}} \AND
\bauthor{\bsnm{Thorbj{\o}rnsen},~\bfnm{Steen}\binits{S.}}
(\byear{2005}).
\btitle{A new application of random matrices}. 
\bjournal{Ann. of Math. (2)}
\bvolume{162}
\bpages{711--775}.
\bid{doi={10.4007/annals.2005.162.711}, issn={0003-486X}, mr={2183281}}
\end{barticle}
\bptok{imsref}%
\endbibitem

\bibitem{hornjohnson85}
\begin{bbook}[mr]
\bauthor{\bsnm{Horn},~\bfnm{Roger~A.}\binits{R.~A.}} \AND
\bauthor{\bsnm{Johnson},~\bfnm{Charles~R.}\binits{C.~R.}}
(\byear{1985}).
\btitle{Matrix Analysis}.
\bpublisher{Cambridge Univ. Press},
\blocation{Cambridge}.
\bid{doi={10.1017/CBO9780511810817}, mr={0832183}}
\end{bbook}
\bptok{imsref}%
\endbibitem

\bibitem{kantorovich48}
\begin{barticle}[mr]
\bauthor{\bsnm{Kantorovi{\v{c}}},~\bfnm{L.~V.}\binits{L.~V.}}
(\byear{1948}).
\btitle{Functional analysis and applied mathematics}.
\bjournal{Uspehi Matem. Nauk (N.S.)}
\bvolume{3}
\bpages{89--185}.
\bid{issn={0042-1316}, mr={0027947}}
\end{barticle}
\bptok{imsref}%
\endbibitem

\bibitem{kargin11b}
\begin{barticle}[mr]
\bauthor{\bsnm{Kargin},~\bfnm{Vladislav}\binits{V.}}
(\byear{2012}).
\btitle{A concentration inequality and a local law for the sum of two random matrices}.
\bjournal{Probab. Theory Related Fields}
\bvolume{154}
\bpages{677--702}.
\bid{doi={10.1007/s00440-011-0381-4}, issn={0178-8051}, mr={3000559}}
\end{barticle}
\bptok{imsref}%
\endbibitem

\bibitem{kargin13a}
\begin{barticle}[mr]
\bauthor{\bsnm{Kargin},~\bfnm{V.}\binits{V.}}
(\byear{2013}).
\btitle{An inequality for the distance between densities of free convolutions}.
\bjournal{Ann. Probab.}
\bvolume{41}
\bpages{3241--3260}.
\bid{doi={10.1214/12-AOP756}, issn={0091-1798}, mr={3127881}}
\bptnote{check year}%
\end{barticle}
\bptok{imsref}%
\endbibitem

\bibitem{knowlesyin12b}
\begin{barticle}[mr]
\bauthor{\bsnm{Knowles},~\bfnm{Antti}\binits{A.}} \AND
\bauthor{\bsnm{Yin},~\bfnm{Jun}\binits{J.}}
(\byear{2014}).
\btitle{The outliers of a deformed Wigner matrix}.
\bjournal{Ann. Probab.}
\bvolume{42}
\bpages{1980--2031}.
\bid{mr={3262497}}
\end{barticle}
\bptok{imsref}%
\endbibitem

\bibitem{knowlesyin12a}
\begin{barticle}[mr]
\bauthor{\bsnm{Knowles},~\bfnm{Antti}\binits{A.}} \AND
\bauthor{\bsnm{Yin},~\bfnm{Jun}\binits{J.}}
(\byear{2013}).
\btitle{The isotropic semicircle law and deformation of {W}igner matrices}.
\bjournal{Comm. Pure Appl. Math.}
\bvolume{66}
\bpages{1663--1750}.
\bid{doi={10.1002/cpa.21450}, issn={0010-3640}, mr={3103909}}
\bptnote{check year}%
\end{barticle}
\bptok{imsref}%
\endbibitem

\bibitem{maida07}
\begin{barticle}[mr]
\bauthor{\bsnm{Ma{\"{\i}}da},~\bfnm{Myl{\`e}ne}\binits{M.}}
(\byear{2007}).
\btitle{Large deviations for the largest eigenvalue of rank one deformations of {G}aussian ensembles}.
\bjournal{Electron. J. Probab.}
\bvolume{12}
\bpages{1131--1150 (electronic)}.
\bid{doi={10.1214/EJP.v12-438}, issn={1083-6489}, mr={2336602}}
\end{barticle}
\bptok{imsref}%
\endbibitem

\bibitem{male11}
\begin{barticle}[mr]
\bauthor{\bsnm{Male},~\bfnm{Camille}\binits{C.}}
(\byear{2012}).
\btitle{The norm of polynomials in large random and deterministic matrices}.
\bjournal{Probab. Theory Related Fields}
\bvolume{154}
\bpages{477--532}.
\bnote{With an appendix by Dimitri Shlyakhtenko}.
\bid{doi={10.1007/s00440-011-0375-2}, issn={0178-8051}, mr={3000553}}
\end{barticle}
\bptok{imsref}%
\endbibitem

\bibitem{nicaspeicher06}
\begin{bbook}[mr]
\bauthor{\bsnm{Nica},~\bfnm{Alexandru}\binits{A.}} \AND
\bauthor{\bsnm{Speicher},~\bfnm{Roland}\binits{R.}}
(\byear{2006}).
\btitle{Lectures on the Combinatorics of Free Probability}.
\bseries{London Mathematical Society Lecture Note Series}
\bvolume{335}.
\bpublisher{Cambridge Univ. Press},
\blocation{Cambridge}.
\bid{doi={10.1017/CBO9780511735127}, mr={2266879}}
\end{bbook}
\bptok{imsref}%
\endbibitem

\bibitem{pasturvasilchuk00}
\begin{barticle}[mr]
\bauthor{\bsnm{Pastur},~\bfnm{L.}\binits{L.}} \AND
\bauthor{\bsnm{Vasilchuk},~\bfnm{V.}\binits{V.}}
(\byear{2000}).
\btitle{On the law of addition of random matrices}.
\bjournal{Comm. Math. Phys.}
\bvolume{214}
\bpages{249--286}.
\bid{doi={10.1007/s002200000264}, issn={0010-3616}, mr={1796022}}
\end{barticle}
\bptok{imsref}%
\endbibitem

\bibitem{peche06}
\begin{barticle}[mr]
\bauthor{\bsnm{P{\'e}ch{\'e}},~\bfnm{S.}\binits{S.}}
(\byear{2006}).
\btitle{The largest eigenvalue of small rank perturbations of {H}ermitian random matrices}.
\bjournal{Probab. Theory Related Fields}
\bvolume{134}
\bpages{127--173}.
\bid{doi={10.1007/s00440-005-0466-z}, issn={0178-8051}, mr={2221787}}
\end{barticle}
\bptok{imsref}%
\endbibitem

\bibitem{peng12}
\begin{bmisc}[auto:STB|2014/02/12|14:17:21]
\bauthor{\bsnm{Peng},~\bfnm{Minyu}\binits{M.}}
(\byear{2012}).
\bhowpublished{Eigenvalues of deformed random matrices.
Available at \arxivurl{arXiv:1205.0572}.}
\end{bmisc}
\bptok{imsref}%
\endbibitem

\bibitem{pizzorenfrewsoshnikov11}
\begin{barticle}[mr]
\bauthor{\bsnm{Pizzo},~\bfnm{Alessandro}\binits{A.}},
\bauthor{\bsnm{Renfrew},~\bfnm{David}\binits{D.}} \AND
\bauthor{\bsnm{Soshnikov},~\bfnm{Alexander}\binits{A.}}
(\byear{2013}).
\btitle{On finite rank deformations of {W}igner matrices}.
\bjournal{Ann. Inst. Henri Poincar\'e Probab. Stat.}
\bvolume{49}
\bpages{64--94}.
\bid{doi={10.1214/11-AIHP459}, issn={0246-0203}, mr={3060148}}
\end{barticle}
\bptok{imsref}%
\endbibitem

\bibitem{renfrewsoshnikov12}
\begin{barticle}[mr]
\bauthor{\bsnm{Renfrew},~\bfnm{David}\binits{D.}} \AND
\bauthor{\bsnm{Soshnikov},~\bfnm{Alexander}\binits{A.}}
(\byear{2013}).
\btitle{On finite rank deformations of {W}igner matrices {II}: {D}elocalized perturbations}.
\bjournal{Random Matrices Theory Appl.}
\bvolume{2}
\bpages{1250015, 36}.
\bid{doi={10.1142/S2010326312500153}, issn={2010-3263}, mr={3039820}}
\end{barticle}
\bptok{imsref}%
\endbibitem

\bibitem{shiryaev96}
\begin{bbook}[mr]
\bauthor{\bsnm{Shiryaev},~\bfnm{A.~N.}\binits{A.~N.}}
(\byear{1996}).
\btitle{Probability},
\bedition{2nd} ed.
\bseries{Graduate Texts in Mathematics}
\bvolume{95}.
\bpublisher{Springer},
\blocation{New York}.
\bid{doi={10.1007/978-1-4757-2539-1}, mr={1368405}}
\end{bbook}
\bptok{imsref}%
\endbibitem

\bibitem{speicher93}
\begin{barticle}[mr]
\bauthor{\bsnm{Speicher},~\bfnm{Roland}\binits{R.}}
(\byear{1993}).
\btitle{Free convolution and the random sum of matrices}.
\bjournal{Publ. Res. Inst. Math. Sci.}
\bvolume{29}
\bpages{731--744}.
\bid{doi={10.2977/prims/1195166573}, issn={0034-5318}, mr={1245015}}
\end{barticle}
\bptok{imsref}%
\endbibitem

\bibitem{tillmann53}
\begin{barticle}[mr]
\bauthor{\bsnm{Tillmann},~\bfnm{Heinz-G{\"u}nther}\binits{H.-G.}}
(\byear{1953}).
\btitle{Randverteilungen analytischer {F}unktionen und {D}istributionen}.
\bjournal{Math. Z.}
\bvolume{59}
\bpages{61--83}.
\bid{issn={0025-5874}, mr={0057345}}
\end{barticle}
\bptok{imsref}%
\endbibitem

\bibitem{voiculescu91}
\begin{barticle}[mr]
\bauthor{\bsnm{Voiculescu},~\bfnm{Dan}\binits{D.}}
(\byear{1991}).
\btitle{Limit laws for random matrices and free products}.
\bjournal{Invent. Math.}
\bvolume{104}
\bpages{201--220}.
\bid{doi={10.1007/BF01245072}, issn={0020-9910}, mr={1094052}}
\end{barticle}
\bptok{imsref}%
\endbibitem

\bibitem{voiculescu93}
\begin{barticle}[mr]
\bauthor{\bsnm{Voiculescu},~\bfnm{Dan}\binits{D.}}
(\byear{1993}).
\btitle{The analogues of entropy and of {F}isher's information measure in free probability theory. {I}}.
\bjournal{Comm. Math. Phys.}
\bvolume{155}
\bpages{71--92}.
\bid{issn={0010-3616}, mr={1228526}}
\end{barticle}
\bptok{imsref}%
\endbibitem

\end{thebibliography}
\end{document}